\documentclass[11pt]{article}

\usepackage{amsmath,amsfonts,amssymb,amsthm,enumerate,graphicx}
\usepackage{yhmath}

\usepackage{tikz}
\usetikzlibrary{calc}
\usepgflibrary{shapes.geometric}
\usepgflibrary{shapes.misc}
\usepgflibrary{shapes.multipart}
\usetikzlibrary{positioning}
\usetikzlibrary{decorations}
\usetikzlibrary{fit}

\newcommand{\N}{\mathbb{N}}
\newcommand{\C}{\mathbb{C}}
\newcommand{\Z}{\mathbb{Z}}
\newcommand{\K}{\mathbb{K}}
\newcommand{\st}{\mid}
\newcommand{\inv}[1]{#1^{-1}}
\newcommand{\til}{\widetilde}
\newcommand{\ep}{\epsilon}
\newcommand{\union}{\cup}
\newcommand{\intersect}{\cap}
\newcommand{\disjunion}{\sqcup}
\newcommand{\isomorphic}{\cong}
\newcommand{\defeq}{:=}
\newcommand{\compose}{\circ}
\newcommand{\cross}{\times}
\newcommand{\curly}[1]{\mathcal{#1}}
\newcommand{\card}[1]{\lvert #1 \rvert}
\newcommand{\lgen}{\ensuremath \mathopen{<}} 
\newcommand{\rgen}{\ensuremath \mathopen{>}}
\newcommand{\map}{\longrightarrow}
\newcommand{\isomap}{\stackrel{\sim}{\map}}
\newcommand{\mapname}[1]{\stackrel{#1}{\map}}

\newcommand{\multdot}{}

\newcommand{\OG}[3]{\mathcal O_{#1}(G( #2 , #3 ))}
\newcommand{\OM}[3]{\mathcal O_{#1}(M(#2 , #3))}
\newcommand{\x}[3]{x^{#1}_{#2 #3}}
\newcommand{\tx}[3]{\hat{x}^{#1}_{#2 #3}}
\newcommand{\Tx}[3]{T(x^{#1}_{#2 #3})}
\newcommand{\M}[1]{M_{#1}}
\newcommand{\locOG}[1]{ \OG{q}{m}{n} [ \inv{ [\M{#1}] } ] }
\newcommand{\dehom}[1]{\K[\x{#1}{i}{j}][y_{#1}^{\pm1};\sigma_{#1}]}
\newcommand{\qminor}[3]{\left[ #1 \left|  \right. #2 \right]_{\x{#3}{}{}}}
\newcommand{\tqminor}[3]{\left[ #1 \left| \right. #2 \right]_{\tx{#3}{}{}}}

\newcommand{\A}[1]{A_{#1}}
\newcommand{\aalpha}{\tilde{\alpha}}
\newcommand{\dialpha}{w_{\alpha,i}}
\newcommand{\djalpha}{z_{\alpha,j}}
\newcommand{\caseif}{\text{if}\hgap}
\DeclareMathOperator{\content}{content}

\newcommand{\hgap}{\;\;\;}

\theoremstyle{definition}
\newtheorem{defn}{Definition}[section]
\newtheorem{remark}[defn]{Remark}

\theoremstyle{plain}
\newtheorem{lem}[defn]{Lemma}
\newtheorem{prop}[defn]{Proposition}
\newtheorem{thm}[defn]{Theorem}
\newtheorem{cor}[defn]{Corollary}

\newcommand{\BDEFN}{\begin{defn}}
\newcommand{\EDEFN}{\end{defn}}
\newcommand{\BLEM}{\begin{lem}}
\newcommand{\ELEM}{\end{lem}}
\newcommand{\BPROP}{\begin{prop}}
\newcommand{\EPROP}{\end{prop}}
\newcommand{\BTHM}{\begin{thm}}
\newcommand{\ETHM}{\end{thm}}
\newcommand{\BCOR}{\begin{cor}}
\newcommand{\ECOR}{\end{cor}}
\newcommand{\BRMK}{\begin{remark}}
\newcommand{\ERMK}{\end{remark}}

\newcommand{\BPF}{\begin{proof}}
\newcommand{\EPF}{\end{proof}}
\newcommand{\nn}{\nonumber}
\newcommand{\BEQ}{\begin{equation}}
\newcommand{\EEQ}{\end{equation}}
\newcommand{\BEQAR}{\begin{eqnarray}}
\newcommand{\EEQAR}{\end{eqnarray}}
\newcommand{\BCASE}{\begin{cases}}
\newcommand{\ECASE}{\end{cases}}

\title{A quantum analogue of the dihedral action on Grassmannians}

\author{Justin M.\ Allman\footnote{Email: jallman@email.unc.edu} \\
		\textit{Department of Mathematics, the University of North Carolina at Chapel Hill} \\
		\textit{CB \#3250, Phillips Hall, Chapel Hill, NC 27599, USA}
		\and Jan E.\ Grabowski\footnote{Email: j.grabowski@lancaster.ac.uk} \\
		\textit{Department of Mathematics and Statistics, Lancaster University} \\
		\textit{Lancaster, LA1 4YF, United Kingdom}}

\date{17th November 2011}

\begin{document}

\maketitle

\begin{abstract}
\noindent In recent work, Launois and Lenagan have shown how to construct a cocycle twisting of the quantum Grassmannian and an isomorphism of the twisted and untwisted algebras that sends a given quantum minor to the minor whose index set is permuted according to the $n$-cycle $c=(1\,2\, \cdots \,n)$, up to a power of $q$.  This twisting is needed because $c$ does not naturally induce an automorphism of the quantum Grassmannian, as it does classically and semi-classically.  

We extend this construction to give a quantum analogue of the action on the Grassmannian of the dihedral subgroup of $S_{n}$ generated by $c$ and $w_{0}$, the longest element, and this analogue takes the form of a groupoid.  We show that there is an induced action of this subgroup on the torus-invariant prime ideals of the quantum Grassmannian and also show that this subgroup acts on the totally nonnegative and totally positive Grassmannians.  Then we see that this dihedral subgroup action exists classically, semi-classically (by Poisson automorphisms and anti-automorphisms, a result of Yakimov) and in the quantum and nonnegative settings.
\end{abstract}

\section{Introduction}

The Grassmannian $G(m,n)$ of $m$-dimensional subspaces of an $n$-dimensional vector space is an important geometric and algebraic object that occurs in many different contexts.  It is a projective variety and its geometric structure is now well-understood, including in particular a cell decomp\-osition.  It is also well-known that the Grassmannian admits an action of the symmetric group $S_{n}$.  Each point in the Grassmannian may be specified by an $m\cross n$ matrix of rank $m$ and the symmetric group action is by permutation of the columns of this matrix.  This action of course extends to the coordinate ring of the Grassmannian, $\curly{O}(G(m,n))$.

Unfortunately, several other important relations of the Grassmannian do not admit a symmetric group action.  Firstly, the (real) totally nonnegative Grassmannian $G^{\mathrm{tnn}}(m,n)$ does not.  Recall that a real $m \cross n$ matrix is totally nonnegative if all of its $m\cross m$ minors are nonnegative.  The totally nonnegative Grassmannian is the space of all of these modulo the action of $\mathrm{GL}^{+}(m)$, the group of real $m\cross m$ matrices with positive determinant.  Clearly the property of being a point in $G^{\mathrm{tnn}}(m,n)$ is not preserved under arbitrary column permutation.  Secondly, in recent work Launois and Lenagan (\cite{LaunoisLenagan}) showed that the natural action of the Coxeter element $c=(1\,2\, \cdots\, n)$ on indexing sets of quantum minors does not induce an automorphism of the quantum Grassmannian $\OG qmn$.

Not all is lost, however.  Postnikov (\cite[Remark~3.3]{Postnikov}) has observed that one may define an action of $c$ on the totally nonnegative Grassmannian by column permutation and a suitable sign correction.  Yakimov (\cite[Theorem~0.1]{Yakimov-Cyclicity}) has shown that $c$ induces a Poisson automorphism of the complex Grassmannian $G(m,n)_{\C}$, taken with its standard Poisson structure.  Indeed, it has long been known that the longest element $w_{0}$ of the symmetric group induces a Poisson anti-automorphism of $G(m,n)_{\C}$ and that $G(m,n)_{\C}$ admits an action of the maximal torus $T_{n}$ of diagonal matrices.  Combining these, Yakimov notes that the group $I_{2}(n) \ltimes T_{n}$ therefore acts on $G(m,n)_{\C}$ by Poisson automorphisms and anti-automorphisms, where $I_{2}(n)=\lgen c,w_{0} \rgen \isomorphic D_{2n}$ is the dihedral group of order $2n$, given its Coxeter group name.

The paper \cite{LaunoisLenagan} of Launois and Lenagan explains how one may construct a replacement for the cycling automorphism in the quantum setting by twisting the quantum Grassmannian.  In this work, we show how to extend this to a quantum analogue of the above $(I_{2}(n) \ltimes T_{n})$-action and deduce that one obtains a dihedral action on the set of torus-invariant prime ideals of $\OG qmn$.  We also note that one may define an $I_{2}(n)$-action on both the totally nonnegative and totally positive Grassmannians, though these do not admit the torus action.

Further reasons for interest in a dihedral action on the quantum Grassmannian come from the study of quasi-commuting sets of quantum minors and cluster algebras.  Leclerc and Zelevinsky (\cite{LeclercZelevinsky}) have shown that two quantum minors in $\OG qmn$ quasi-commute if and only if the column sets defining them satisfy a combinatorial condition called weak separability.  As noted by Scott (\cite[Proposition~4]{Scott-QMinors}), the natural action of the group $I_{2}(n)$ on $m$-subsets of $\{1,\ldots ,n\}$ preserves weak separability and so the question of an analogue of the dihedral action on $\OG qmn$ naturally arises here too.  

In \cite{Gr2nSchubertQCA}, the second author and Launois have observed the quantum cycling of \cite{LaunoisLenagan} playing a role in a quantum cluster algebra structure on $\OG q3m$ for $m=6,7,8$ and also hints of a quantum dihedral action.  Since quantum clusters are by definition quasi-commuting sets, this is not so surprising, although it should be noted that in the cases mentioned not all quantum cluster variables are quantum minors but the dihedral action is still evident.  Also, Assem, Schiffler and Shramchenko (\cite{AsScSh-clusterautom}) have studied automorphism groups of (unquantized) cluster algebras and shown that a cluster algebra of type $A_{n-3}$ has a dihedral cluster automorphism group of order $2n$.  Fomin and Zelevinsky (\cite{FZ-CA2}) showed that $\OG{}{2}{n}$ is a cluster algebra of type $A_{n-3}$ and hence $\OG{}{2}{n}$ has cluster automorphism group isomorphic to $I_{2}(n)$.  It is expected that the results presented here will aid the understanding of the quantum cluster algebra structures conjectured to exist for all quantum Grassmannians.

As noted above, in the commutative setting the Grassmannian admits a symmetric group action and we are clearly a considerable distance from having a quantum analogue of this, if indeed one exists.  A direction for future work would be to try to extend our dihedral action further.  It would also be interesting to know whether other related geometric results can be similarly improved, such as whether the action of $c$ on the Lusztig strata of the Grassmannian (\cite{KnutsonLamSpeyer},\cite{Yakimov-Cyclicity}) extends to a dihedral action on the strata.  We do not address these questions here, though.

\vfill \pagebreak

The main result of this paper, the quantum analogue of the dihedral action on $\OG qmn$, can be expressed in categorical language.  A groupoid is a category in which every morphism has an inverse and the natural way to regard an algebra and its automorphism group is as a groupoid with one object; the latter is usually simply called a group, of course.  Then what we see is that under quantization the subgroup $I_{2}(n)$ of $\mathrm{Aut}(\OG{}{m}{n})$ is replaced by a groupoid with infinitely many objects, the twists of $\OG qmn$, and arrows generated by the isomorphisms we establish.  This groupoid is equivalent to but not equal to a group; in the quantum setting, we see many algebras isomorphic to $\OG qmn$ but not equal to it, but under passing to the classical limit these algebras coincide once more and we recover $\OG{}{m}{n}$ and this dihedral part of its automorphism group.

We also note that in the semi-classical and quantum settings, it is not enough to consider just automorphisms of algebras.  We need to include anti-automorphisms, such as that which is the quantum analogue of the action of the element $w_{0}$, as described in this work.  This is natural when working with noncommutative structures having a classical limit.  Thus we use the term ``automorphism groupoid'' to encompass automorphisms and anti-automorphisms, for example.

The organisation of this paper is as follows.  In Sections~\ref{ss:defns} and \ref{ss:twisting}, we recall some required definitions.   We then recall the isomorphism of \cite[Theorem~5.9]{LaunoisLenagan} and show that by twisting the algebras related by that map, we can produce a family of algebras, all isomorphic to the quantum Grassmannian, together with a collection of maps between them that have the effect on quantum minors of acting on their indexing sets by the cycle $c$, up to a power of $q$.  That is, in Section~\ref{ss:groupoid.1} we construct the part of our automorphism groupoid consisting of the quantum rotations.  Then in Section~\ref{ss:dihedral} we show that there is an anti-automorphism of the quantum Grassmannian corresponding to the action of $w_{0}$ and use this to complete the construction of our dihedral automorphism groupoid.  In Section~\ref{s:twistingHprimes} we discuss the consequences of these results for the torus-invariant prime ideals of $\OG qmn$.  

In fact, there is a choice involved in the construction we have given, via the isomorphism of Launois and Lenagan.  Specifically, their proof of the existence of this map relies on the dehomogenisation isomorphism of \cite{KLR}, linking a localisation of the quantum Grassmannian $\OG qmn$ at a consecutive minor with a skew-Laurent extension of the quantum matrix algebra $\OM qm{n-m}$.  The dehomogenisation isomorphism exists for any choice of consecutive minor and Launois and Lenagan chose to work with the one with indexing set $\{1,\ldots ,m\}$.  In Section~\ref{s:groupoid.any.integer} we show that the method of Launois and Lenagan, and hence also our own groupoid construction, can be carried out for any choice of consecutive minor.  Hence we produce an integer-parameterized family of dihedral groupoids, where the one in Section~\ref{ss:dihedral} is associated to the parameter being equal to $1$.

Finally, we conclude in Section~\ref{s:tnn} with some comments on the dihedral action on the totally nonnegative and totally positive Grassmannians.

\vspace{-1em}
\subsection*{Acknowledgements}

We would like to thank St\'{e}phane Launois for explanation of the results of \cite{LaunoisLenagan} and Tom Lenagan for several helpful suggestions.

The first named author would like to thank the University of North Carolina at Chapel Hill for a gift of travel funds, as well as Green--Templeton College and the Mathematical Institute of the University of Oxford for provision of facilities. The second author would like to acknowledge the provision of facilities by Keble College and the Mathematical Institute of the University of Oxford, where the majority of this work was carried out.

\section{A dihedral action on quantum Grassmannians}\label{s:quantgrass}

\subsection{Definitions and other prerequisite results}\label{ss:defns}

Throughout we work over a base field, which we will denote by $\K$. We set $\K^\ast \defeq \K\setminus\{0\}$ and let $q\in\K^\ast$.

Let $m , n\in\N$ be such that $m<n$ and assume that there exists an element $p\in\K^\ast$ such that $p^m=q^2$ and $p\neq 1$. The major object under consideration is the so-called \emph{quantum Grassmannian}, which we denote by $\OG qmn$. More precisely, this is the quantized coordinate ring of the space of $m$-dimensional linear subspaces of $\K^n$, which is constructed as follows. 

We recall that the quantum matrix algebra $\OM qmn$, the quantization of the coordinate ring of the affine variety of $m\cross n$ matrices with entries in $\K$, is the $\K$-algebra generated by the set $\{ X_{ij} \mid 1\leq i\leq m,\ 1\leq j \leq n \}$ subject to the quantum $2\cross 2$ matrix relations on each $2\cross 2$ submatrix of \[ \begin{pmatrix} X_{11} & X_{12} & \cdots & X_{1n} \\ \vdots & \vdots & \ddots & \vdots \\ X_{m1} & X_{m2} & \cdots & X_{mn} \end{pmatrix}, \] where the quantum $2\cross 2$ matrix relations on $\left( \begin{smallmatrix} a & b \\ c & d \end{smallmatrix} \right)$ are
\begin{align*} ab & = qba & ac & = qca & bc & = cb \\ bd & = qdb & cd & = qdc & ad-da & = (q-q^{-1})bc. \end{align*}

The quantum Grassmannian $\OG qmn$ is defined to be the subalgebra of the quantum matrix algebra $\OM qmn$ generated by the quantum Pl\"{u}cker coordinates, these being the $m\cross m$ quantum minors of $\OM qmn$ defined as follows.  The $m\cross m$ quantum minor $\Delta_{q}^{I}$ associated to the $m$-subset $I=\{ i_{1} < i_{2} < \cdots < i_{m} \}$ of $\{1, \dots , n\}$ is defined to be
\[ \Delta_{q}^{I} \defeq \sum_{\sigma \in S_{m}} (-q)^{l(\sigma)}X_{1i_{\sigma(1)}}\cdots X_{m\mspace{0.5mu}i_{\sigma(m)}} \] where $S_{m}$ is the symmetric group of degree $m$ and $l$ is the usual length function on this.  (In fact we are considering the quantum minor $\Delta_{\{1,\ldots,m\}}^{I}$ but since we are working in $\OM qmn$ there is no other choice for the row subset and so we omit it.)  Then we denote by \[ \curly{P}_{q} \defeq \{ \Delta_{q}^{I}  \mid I \subseteq \{1,\ldots ,n\}, \card{I}=m \} \] the set of all quantum Pl\"{u}cker coordinates and this is the defining generating set of $\OG qmn$.

Our notation mostly coincides with that of \cite{LaunoisLenagan}, and in particular we also observe the following notational conveniences. Throughout the sequel, whenever $j\in\Z$, let $\til j$ denote the element of the set $\{1,\ldots, n\}$ which is congruent to $j$ modulo $n$.  Given an $m$-subset $I$ of $\{ 1,\ldots ,n\}$, we denote by $[I]$ the minor $\Delta_{q}^{I}$.  We will also abbreviate operations on indexing sets of minors in the natural way, writing for example $I+a$ (or $a+I$) for $\{ \widetilde{i+a} \mid i\in I \}$.

We say that two elements $x$ and $y$ of $\OG qmn$ quasi-commute if there exists an integer $r$ such that $yx=q^{r}xy$.  When calculating quasi-commutation relations between minors, we will on occasion call upon the following result, a quantum version of Muir's Law of Extensible Minors.  We use the same form of this result as \cite{LaunoisLenagan}, which was first obtained by Krob and Leclerc (\cite[Theorem~3.4]{KrobLeclerc}).

\BPROP \label{prop:QuantumMuirsLaw} Let $I_{s}$, $J_{s}$, for $1\leq s \leq d$, be $m$-element subsets of $\{ 1,\ldots ,n\}$ and let $c_{s}\in \K$ be such that $\sum_{s=1}^{d} c_{s}[I_{s}][J_{s}]=0$ in $\OG qmn$.  Suppose that $P$ is a subset of $\{ 1,\ldots ,n\}$ such that $(\bigcup_{s=1}^{d} I_{s}) \union (\bigcup_{s=1}^{d} J_{s}) \subseteq P$ and let $\bar{P}$ denote $\{ 1,\ldots ,n\}\setminus P$.  Then \[ \sum_{s=1}^{d} c_{s}[I_{s} \sqcup \bar{P}][J_{s} \sqcup \bar{P}]=0 \] holds in $\OG{q}{m'}{n}$, where $m'=m+\card{\bar{P}}$. \qedhere
\EPROP

A consequence of this result is that one may simplify relations by deleting common members of the index sets of the minors involved.

\subsection{Twisting by a 2-cocycle}\label{ss:twisting}

When a $\K$-algebra $A$ is graded by a semigroup, one can twist the multiplication in $A$ by using a 2-cocycle to produce a new multiplication.  In the situation relevant to us, that of $\Z^{n}$-graded algebras, a 2-cocycle with values in $\K^{\ast}$ is a map $\chi\colon \Z^{n}\cross \Z^{n} \to \K^{\ast}$ such that \[ \chi(s,t+u)\chi(t,u)=\chi(s,t)\chi(s+t,u) \] for all $s,t,u\in \Z^{n}$.  Then given a $\Z^{n}$-graded $\K$-algebra $A$ and 2-cocycle $\chi$, one may form the twist $T(A)$ of $A$ by $\chi$ by taking $T(A)$ to be isomorphic to $A$ as a vector space by an isomorphism $a\mapsto T(a)$ and defining the multiplication in $T(A)$ by $T(a)T(b)=\chi(s,t)T(ab)$ for homogeneous elements $a,b \in A$ with multi-degrees $s$ and $t$, respectively.  The cocycle condition above ensures the associativity of the product in the $\K$-algebra $T(A)$. The next lemma will be used a number of times and records the fact that twisting preserves graded isomorphisms.

\BLEM	\label{lemma:twist.preserves.isomorphism}
Let $A$ and $B$ be $\Z^n$-graded $\K$-algebras and $\varphi\colon A \map B$ a graded isomorphism. Then if $A$ and $B$ are both twisted by the same cocycle $\chi$, $T(A)\isomorphic T(B)$ via a map $\widehat\varphi$ given by $T(a)\mapsto T(\varphi(a))$.\qed
\ELEM

The quantum Grassmannian $\OG qmn$ has a natural $\Z^n$-grading given by
\BEQ\label{eqn:natural.grading}
\content([I])=\sum_{i\in I} \ep(i)
\EEQ
where $[I]$ is a generating quantum minor and $\{\ep(1),\ldots,\ep(n)\}$ denotes the standard basis of $\Z^n$. As shown in \cite[Lemma~5.3]{LaunoisLenagan}, the map $\Gamma\colon \Z^n\times\Z^n\map\K^\ast$ given by
	\BEQ	
	\Gamma((s_1,\ldots,s_n),(t_1,\ldots,t_n))=\prod_{j\neq n} p^{s_n t_j}	\nn
	\EEQ
is a 2-cocycle. Throughout the sequel, the twist of an algebra $A$ by $\Gamma$ will be denoted $T(A)$. Additionally, we will need the cocycle $\gamma\colon\Z^n\times\Z^n\map\K^\ast$ given by
	\BEQ
	\gamma((s_1,\ldots,s_n),(t_1,\ldots,t_n))=\prod_{j\neq 1} (1/p)^{s_1 t_j}.\nn
	\EEQ
One can check that this is a cocycle in the same way as $\Gamma$. Henceforth, the twist of $A$ by $\gamma$ will be denoted $\tau(A)$.  We will also write $\Gamma(a,b)$ (respectively $\gamma(a,b)$) for $\Gamma(\content(a),\content(b))$ (resp.\ $\gamma(\content(a),\content(b))$) for homogeneous elements $a$ and $b$.

On the level of elements, the notation $T(\ )$ (respectively $\tau(\ )$) can be somewhat cumbersome, so for $a\in A$ we set $\hat{a}=T(a)$ (resp.\ $\hat{a}=\tau(a)$), but retain the former notation when this is clearer.

\subsection{Quantum rotations}\label{ss:groupoid.1}

Our goal in this section is to construct algebras and maps corresponding to powers of the cycle $c$, that is, ``quantum rotations''.  More precisely, using the cocycle $\Gamma$, Launois and Lenagan prove the following theorem. 

\BTHM[{\cite[Theorem 5.9]{LaunoisLenagan}}]\label{thm:LL.theta}
There is an isomorphism $$\theta\colon T(\OG qmn) \isomap \OG qmn$$ such that $T([i_1,\ldots,i_m])$ is sent to the minor $[i_1+1,\ldots,i_m+1]$ when $i_m<n$ and $T([i_1,\ldots,i_{m-1},n])$ is sent to $q^{-2}[1,i_1+1,\ldots,i_{m-1}+1]$. \qed
\ETHM

Thus, on the level of quantum minors, the map $\theta$ acts like the cycle $c$, up to a power of $q$. From this map, we will now construct a family of twists of $\OG qmn$ and algebra isomorphisms which have the effect on minors of permuting the column indices by $c^\ell$ for any integer $\ell$.

We wish to apply Lemma~\ref{lemma:twist.preserves.isomorphism} iteratively to the isomorphism $\theta$. We have described the natural $\Z^n$-grading of $\OG qmn$, so in order to ensure that $\theta$ preserves this grading, we are forced to assign for $T([I])\in T(\OG qmn)$ that $$\content(T([I])) = \content([{I+1}])$$ where the left-hand side describes the grading on $T(\OG qmn)$ and the right-hand side is the grading on $\OG qmn$ previously defined by equation (\ref{eqn:natural.grading}).

\BLEM\label{lem:gamma.Gamma.1}
Let $[I]$ and $[J]$ be any two generating quantum minors of $\OG qmn$. Then $$\gamma([{I+1}],[{J+1}]) \multdot \Gamma([I],[J]) \equiv 1.$$
\ELEM

\BPF
Notice that $\Gamma([I],[J])=1$ unless $n\in I$. Similarly, $\gamma([{I+1}],[{J+1}])=1$ unless $1\in{I+1}$, that is, unless $n\in I$. So when $n\notin I$, both factors are already $1$. 

When $n\in I$, set $e_J = \card{J \setminus \{n\}}$, so that $\Gamma([I],[J])=p^{e_J}$. Note that $e_J$ is also equal to the number of entries in $J+1$ which are distinct from $1$, and therefore $\gamma([{I+1}],[{J+1}]) = (1/p)^{e_J}$.
\EPF

The above lemma implies that the twist of $T(\OG qmn)$ by the cocycle $\gamma$ gives an equality of algebras $\tau(T(\OG qmn)) = \OG qmn$. Using this, and applying Lemma~\ref{lemma:twist.preserves.isomorphism} to the isomorphism in Theorem \ref{thm:LL.theta}, we have proved

\BTHM\label{thm:Theta.1}
Let $I=\{i_1<\cdots<i_m\}\subset\{1,\ldots,n\}$ and set $\lambda_I=1$ if $i_m<n$ and $\lambda_I=q^{-2}$ if $i_m = n$. The map $\OG qmn \map \tau(\OG qmn)$ which sends the generating quantum minor $[I]$ to $\lambda_I \multdot \tau([{I+1}])$ defines an algebra isomorphism $$\Theta_{1}\colon \OG qmn \isomap \tau(\OG qmn).$$ Also, $\Theta_1$ induces a $\Z^{n}$-grading on $\tau(\OG qmn)$ given by $\content(\tau([I]))\defeq \content([{I-1}])$. \qed
\ETHM

Now, let $r$ be a positive integer, and $A$ an algebra graded by $\Z^n$. Let $\tau^r(A)$ denote the algebra $\tau(\tau\cdots(\tau(A))\cdots)$ where $\tau$ appears $r$ times. Similarly, for $a\in A$ let $\tau^r(a)$ denote the corresponding element in $\tau^r(A)$. As is natural, we let $\tau^0(A) = A$. We obtain the following corollaries to Theorem~\ref{thm:Theta.1}.

\BCOR\label{cor:Theta.ell}
Let $\ell$ be a positive integer and $\lambda_I$ defined as above. Then there is a graded algebra isomorphism $$\Theta_\ell:\tau^{\ell-1}(\OG qmn) \isomap \tau^\ell(\OG qmn)$$ which sends the homogeneous generator $\tau^{\ell-1}([I])$ to $\lambda_I\multdot\tau^\ell([{I+1}])$.
\ECOR

\BPF
Induct on $\ell$. For $\ell=1$, we just have the map $\Theta_1$ of Theorem \ref{thm:Theta.1}. For $\ell>1$, the grading on $\tau^{\ell-1}(\OG qmn)$ is assigned by the isomorphism $\Theta_{\ell-1}$ and thus one can twist by $\tau$ once more to obtain $\tau^\ell(\OG qmn)$ so that the target of $\Theta_\ell$ is well defined. The map $\Theta_\ell$ is just given by applying Lemma~\ref{lemma:twist.preserves.isomorphism} to $\Theta_{\ell-1}$. Finally, the resulting isomorphism preserves gradings subject to the assignment $\content(\tau^\ell([I]))\defeq\content(\tau^{\ell-1}([{I-1}]))$.
\EPF

Given a generating quantum minor $[I]\in\OG qmn$, and $r\in\Z_{>0}$ define the scalar $$\Lambda_I(r)\defeq \prod_{s=1}^r \lambda_{c^{s-1}(I)}$$ which is evidently a power of $q$. Finally, we can obtain an isomorphism which, on minors, has the effect of applying $c^r$ to the column indices, up to a power of $q$.

\BCOR\label{cor:Theta.ell.composition}
There is an isomorphism $$\OG qmn \isomap \tau^{\ell}(\OG qmn)$$ given by the composition $\Theta_\ell\Theta_{\ell-1}\cdots\Theta_{1}$, sending the quantum minor $[I]$ to $\Lambda_I(\ell)\multdot\tau^\ell([{I+\ell}])$, or equivalently to $\Lambda_I(\ell)\multdot\tau^\ell([c^\ell(I)])$.\qed
\ECOR

All of the above results were obtained by iteratively applying Lemma~\ref{lemma:twist.preserves.isomorphism} to the isomorphism $\theta$ in Theorem \ref{thm:LL.theta}, twisting both sides by the cocycle $\gamma$. Of course one can play the exact same game using the cocycle $\Gamma$. Doing so, one obtains the following results.

\BCOR\label{cor:Theta.-ell}
For every non-negative integer $\ell$, there is a graded isomorphism $$\Theta_{-\ell}\colon T^{\ell+1}(\OG qmn) \isomap T^{\ell}(\OG qmn)$$ which sends $T^{\ell+1}([I])$ to $\lambda_I\multdot T^{\ell}([{I+1}])$.\qed
\ECOR

Above, the gradings on the source and target are assigned iteratively as in Corollary \ref{cor:Theta.ell}. Explicitly, $$\content(T^\ell([I])) = \content(T^{\ell-1}([{I+1}])) = \cdots = \content([{I+\ell}]).$$

\BRMK
Note that when $\ell=0$ in Corollary \ref{cor:Theta.-ell}, we have exactly the map $\theta$ of Theorem \ref{thm:LL.theta}. That is, $\Theta_0$ is the isomorphism $\theta$ of Launois and Lenagan on which the entire construction is based.
\ERMK

\BCOR\label{cor:Theta.-ell.composition}
There is a graded isomorphism $$T^{\ell}(\OG qmn) \isomap \OG qmn$$ given by the composition $\Theta_{0}\Theta_{-1}\cdots\Theta_{-\ell+1}$, sending $T^{\ell}([I])$ to $\Lambda_I(\ell)\multdot[{I+\ell}]=\Lambda_I(\ell)\multdot[c^\ell(I)]$.\qed
\ECOR

Now we have maps $\Theta_{\ell}$ for each $\ell \in \Z$, each of which corresponds to the action of $c$ on minors, and the composition of $r$ consecutive such maps corresponds to the action of $c^{r}$.  That is, we have produced the part of the groupoid which corresponds to the cyclic subgroup $\lgen c \rgen$ of $I_{2}(n)$.

It is natural to inquire about the quantum analogue of the relation $c^n=1$. The answer is given by the following proposition.

\BPROP
The scalar $\Lambda_I(n)$ is independent of $I$. Explicitly, $\Lambda_I(n) = q^{-2m}$ for any $I$.
\EPROP

\BPF
Recall that $\lambda_J$ takes the value $q^{-2}$ if $n\in J$, otherwise it is $1$. Consider an $m$-subset $I=\{i_1,\ldots,i_m\}\subset\{1,\ldots,n\}$. Thus, to compute $\Lambda_I(n) = \prod_{s=1}^n \lambda_{c^{s-1}(I)}$, we must count the number of times that $n$ appears in $c^{r}(I)$ as $r$ ranges from $0$ to $n-1$. This is easily seen to be exactly $m$ times, independent of $I$, since each $i_k\in I$ is cycled to the value $n$ exactly once.
\EPF

Consequently, for $\ell=n$ the graded isomorphisms described in Corollary~\ref{cor:Theta.ell.composition} and Corollary~\ref{cor:Theta.-ell.composition} scale every quantum minor by the same scalar.  These maps are not themselves scalar multiples of the identity but their classical limits are exactly the identity map on $\OG{}{m}{n}$.  Indeed, it is easy to see that any composition of the form $\Theta_{n+r-1}\Theta_{n+r-2}\cdots \Theta_{r+1}\Theta_{r}$, that is $n$ consecutive $\Theta_{i}$'s, also has this property, as such a composition also takes a quantum minor to itself multiplied by the scalar $\Lambda_{I}(n)$.

\subsection{Quantum reflections and dihedral relations}\label{ss:dihedral}

The family of graded algebras $\{\OG qmn, T^i(\OG qmn), \tau^j(\OG qmn) \colon i,j\in\Z_{>0}\}$ will form the complete set of objects in the groupoid we are constructing. Together, the arrows defined by the isomorphisms $\{\Theta_\ell \colon \ell\in\Z\}$ should be viewed as ``quantum rotations" and within the groupoid they are the quantum analogue of the action of the cyclic subgroup in $I_2(n)$ generated by $c$. It remains to construct arrows which will serve as the ``quantum reflections", that is, arrows which on minors have the effect of permuting the column indices by the longest element $w_{0}\in S_n$.

Recall that if $i\in\{1,\ldots,n\}$ then $w_{0}(i) = n - i + 1$.  We construct the quantum analogue of the action of this element from two maps between $\OG qmn$ and $\OG {\inv q}mn$.  The first comes from a well-known anti-isomorphism $f \colon  \OM qmn \to \OM {\inv q}mn$ sending $X_{ij}$ to $Y_{ij}$, where we take generators $Y_{ij}$ for $\OM {\inv q}mn$ in exactly the same manner as for $\OM qmn$.  As noted by Lenagan and Rigal in \cite[Proposition~3.4.1]{LR-qGrAlgs}, this map has the property that $f([I])=[I]$ for every minor $I$ and so induces an anti-isomorphism $f' \colon \OG qmn \to \OG {\inv q}mn$ with the same property.  (We will not distinguish by notation which algebra a minor belongs to, as the context will always make this clear.)

Next, Kelly, Lenagan and Rigal remark in \cite{KLR} that the element $w_{0}$ induces an isomorphism of quantum Grassmannians \[ g \colon \OG qmn \isomap \OG{\inv q}mn \] (their map is called $\delta$) sending the generating minor $[I]$ to the minor $[w_{0}(I)]$.  Hence the composition $\Omega_{0} \defeq \inv{(f')} \circ g$ is an anti-automorphism of $\OG qmn$ with $\Omega_{0}([I])=[w_{0}(I)]$.

Our next goal is to show that $w_{0}$ also provides a correspondence between the higher twisted powers $T^{\ell}(\OG qmn)$ and $\tau^{\ell}(\OG qmn)$ for each positive integer $\ell$.

\BPROP\label{prop:Omega.ell}
Let $\ell$ be a positive integer. There exists an anti-isomorphism $$\Omega_\ell\colon T^\ell(\OG qmn) \isomap \tau^{\ell}(\OG qmn)$$ given by $$T^{\ell}([I]) \mapsto \Lambda_I(\ell) \multdot \Lambda_{(w_{0} c^\ell) I}(\ell) \multdot \tau^\ell([w_{0}(I)]).$$
\EPROP

\BPF
Define the map $\Omega_\ell$ as the composition $$(\Theta_{\ell} \cdots \Theta_{2} \Theta_1)\Omega_0(\Theta_{0} \Theta_{-1} \cdots \Theta_{-\ell+1}).$$
Tracing $T^{\ell}([I])$ through the sequence above with Corollaries \ref{cor:Theta.ell.composition} and \ref{cor:Theta.-ell.composition} plus the dihedral group relation $c^\ell w_{0} c^\ell = w_{0}$ yields the result. 
\EPF

We may picture our quantum dihedral groupoid as follows:

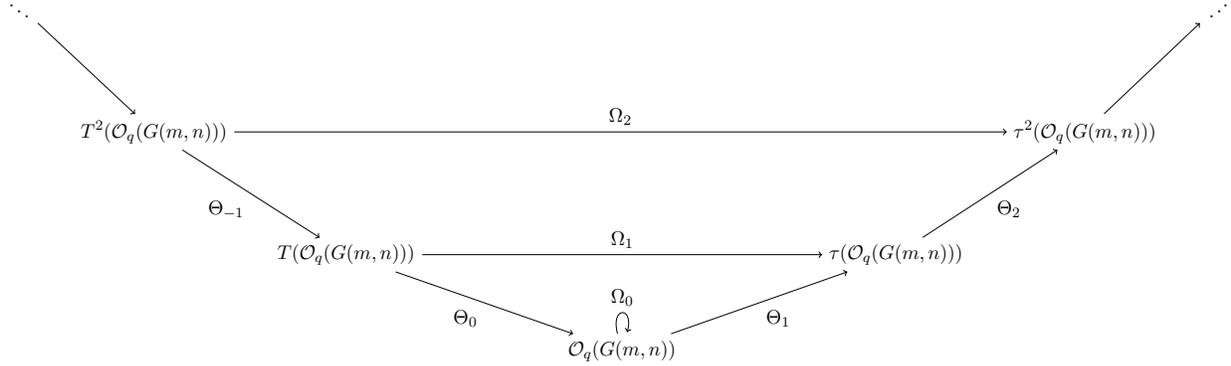
\begin{center} 
\scalebox{0.682}{\begin{tikzpicture}[auto,node distance=1cm]

\node (T0) {$\OG{q}{m}{n}$};
\node (T1) [above left=of T0,xshift=-2cm,yshift=0.5cm] {$T(\OG{q}{m}{n})$};
\node (tau1) [above right=of T0,xshift=2cm,yshift=0.5cm] {$\tau(\OG{q}{m}{n})$};
\node (T2) [above left=of T1,yshift=1cm] {$T^{2}(\OG{q}{m}{n})$};
\node (T3) [above left=of T2,yshift=1cm] {$\ddots$};
\node (tau2) [above right=of tau1,yshift=1cm] {$\tau^{2}(\OG{q}{m}{n})$};
\node (tau3) [above right=of tau2,yshift=1cm] {$\adots$};

\draw[loop above] (T0) to node {$\Omega_{0}$} (T0);

\draw[->] (tau1) to node[swap] {$\Theta_{2}$} (tau2);
\draw[->] (T0) to node[swap] {$\Theta_{1}$} (tau1);
\draw[->] (T1) to node[swap] {$\Theta_{0}$} (T0);
\draw[->] (T2) to node[swap] {$\Theta_{-1}$} (T1);
\draw[->] (T3) to (T2);
\draw[->] (tau2) to (tau3);

\path[->,relative] (T2) edge node {$\Omega_{2}$} (tau2);
\path[->,relative] (T1) edge node {$\Omega_{1}$} (tau1);

\end{tikzpicture}}\label{groupoid}
\end{center}

\noindent We have included in this picture only the maps described above; of course, many more maps are defined within the groupoid as compositions of these and their inverses.  In particular, the identity map at each object is included.  In the following section, we demonstrate the consequences of these isomorphisms for the prime spectra of quantum Grassmannians.

\section{A dihedral action on the $\curly{H}$-prime spectrum}\label{s:twistingHprimes}

Assume that $q$ is not a root of unity.  Then Launois, Lenagan and Rigal (\cite{LLR}) have shown that the prime ideals of $\OG qmn$ are completely prime and furthermore the prime spectrum of $\OG qmn$ has a stratification, in the sense of Goodearl and Letzter (see for example \cite{Brown-Goodearl}), parameterized by a special class of prime ideals of $\OG qmn$, the $\curly{H}$-prime ideals.  Letting $\curly{H}=(\K^{\ast})^{n}$ be an algebraic torus, the $\Z^{n}$-grading on $\OG qmn$ induces a rational action of $\curly{H}$ by \[ (h_{1},\ldots ,h_{n})\cdot [\{ i_{1},\ldots ,i_{m}\}]=h_{i_{1}}h_{i_{2}}\cdots h_{i_{m}}[\{i_{1},\ldots ,i_{m}\}]. \]
The homogeneous prime ideals of $\OG qmn$ are exactly the primes invariant under this action.  These distinguished primes are called $\curly{H}$-primes and the set of all $\curly{H}$-primes, $\curly{H}\text{-Spec}(\OG qmn)$, is called the $\curly{H}$-prime spectrum of $\OG qmn$.  Importantly, this set has been shown in \cite{LLR} to be finite: excluding the augmentation ideal, the $\curly{H}$-primes are in bijection with a set of Cauchon diagrams on certain Young diagrams (see \cite[Section~2]{LLR} for definitions and examples).  The number of the latter and hence the cardinality of $\curly{H}\text{-Spec}(\OG qmn)$ is known, following work of Williams (\cite{Williams}).

Additionally, Launois and Lenagan (\cite[Section~6]{LaunoisLenagan}) have shown that if $P$ is an $\curly{H}$-prime ideal of $\OG qmn$ then $T(P)$ is an $\curly{H}$-prime ideal of $T(\OG qmn)$ and that the isomorphism $\Theta_{0}$ (their $\theta$) yields a self-bijection of the $\curly{H}$-prime spectrum of $\OG qmn$.  Furthermore, for $P$ an $\curly{H}$-prime ideal, $[I]\in P$ if and only if $[I+1]\in \Theta_{0}(T(P))$ so the sets of quantum minors that are in $\curly{H}$-prime ideals are permuted by twisting and $\Theta_{0}$.  By repeating this, we see immediately that if $[I]\in P$ then $[I+\alpha]$ lies in some $\curly{H}$-prime, for any $\alpha \in \Z$.  Our Corollaries~\ref{cor:Theta.ell.composition} and \ref{cor:Theta.-ell.composition} provide an alternative proof of this, by identical arguments to the results of Launois and Lenagan.

Since we have a genuine anti-automorphism $\Omega_{0}\colon \OG qmn \to \OG qmn$, it is immediate that $\Omega_{0}$ induces a self-bijection of the set of $\curly{H}$-primes of $\OG qmn$.  Then, as $\Omega_{0}([I])=[w_{0}(I)]$ for all minors, the sets of quantum minors in $\curly{H}$-prime ideals are permuted according to the action of $w_{0}$ on their indexing sets.  Thus we see that we have an action of the group $I_{2}(n)$ on $\curly{H}\text{-Spec}(\OG qmn)$ and we expect this to be useful in furthering understanding of this set.

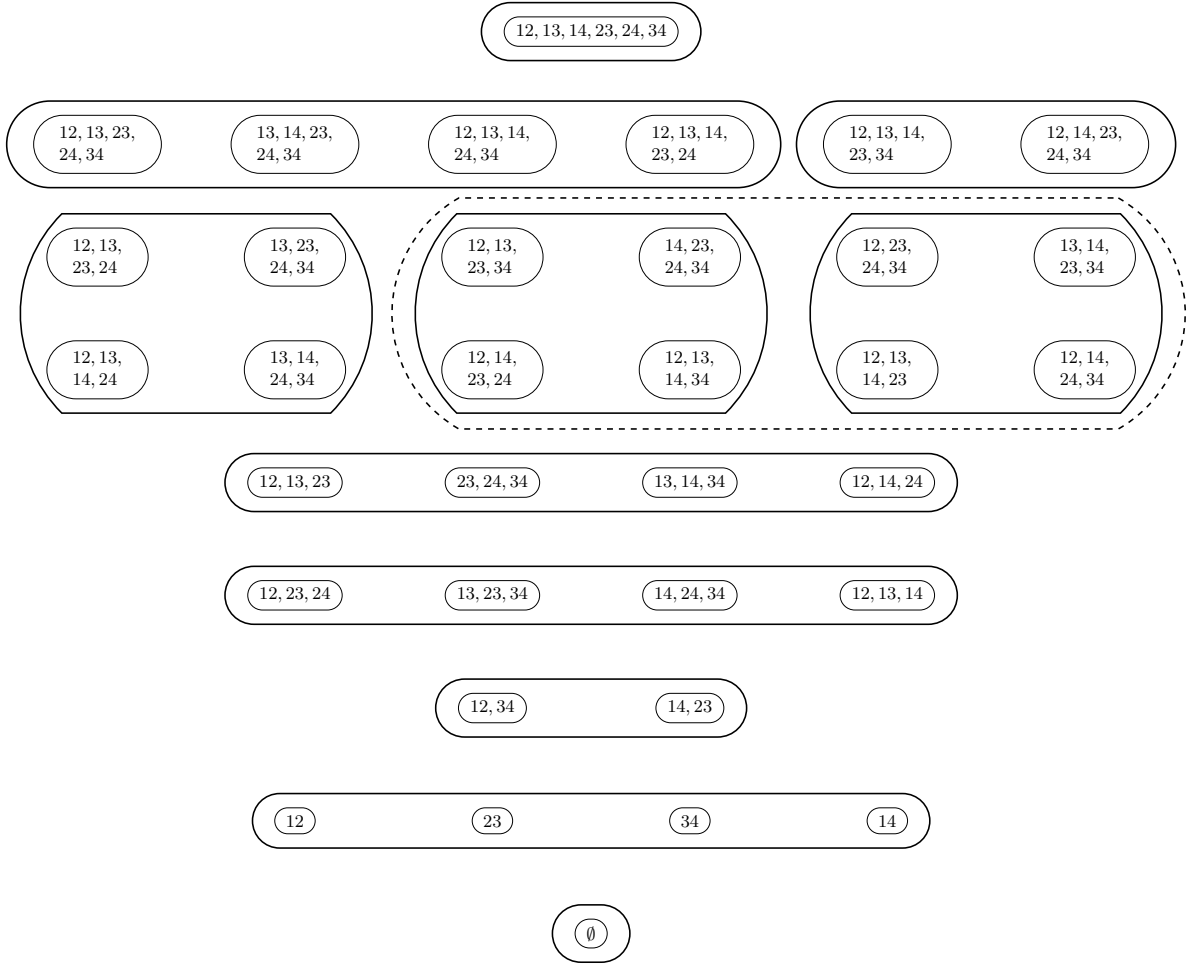
\begin{figure}[ht]
\begin{center} 
\scalebox{0.75}{\begin{tikzpicture}[x=3.5cm,y=2cm,font=\footnotesize]

\node (61) [rounded rectangle,draw] at (2.5,8) {$12,13,14,23,24,34$};
\node[inner sep=7pt,draw,thick,rounded rectangle,fit=(61)] {};

\node (51) [rounded rectangle,draw] at (0,7) {$\begin{array}{l} 12,13,23, \\ 24,34 \end{array}$};
\node (52) [rounded rectangle,draw] at (1,7) {$\begin{array}{l} 13,14,23, \\ 24,34 \end{array}$};
\node (53) [rounded rectangle,draw] at (2,7) {$\begin{array}{l} 12,13,14, \\ 24,34 \end{array}$};
\node (54) [rounded rectangle,draw] at (3,7) {$\begin{array}{l} 12,13,14, \\ 23,24 \end{array}$};
\node[inner sep=7pt,draw,thick,rounded rectangle,fit=(51) (52) (53) (54)] {};
\node (55) [rounded rectangle,draw] at (4,7) {$\begin{array}{l} 12,13,14, \\ 23,34 \end{array}$};
\node (56) [rounded rectangle,draw] at (5,7) {$\begin{array}{l} 12,14,23, \\ 24,34 \end{array}$};
\node[inner sep=7pt,draw,thick,rounded rectangle,fit=(55) (56)] {};

\node (41) [rounded rectangle,draw] at (0,6) {$\begin{array}{l} 12,13, \\ 23,24 \end{array}$};
\node (42) [rounded rectangle,draw] at (1,6) {$\begin{array}{l} 13,23, \\ 24,34 \end{array}$};
\node (43) [rounded rectangle,draw] at (2,6) {$\begin{array}{l} 12,13, \\ 23,34 \end{array}$};
\node (44) [rounded rectangle,draw] at (3,6) {$\begin{array}{l} 14,23, \\ 24,34 \end{array}$};
\node (45) [rounded rectangle,draw] at (4,6) {$\begin{array}{l} 12,23, \\ 24,34 \end{array}$};
\node (46) [rounded rectangle,draw] at (5,6) {$\begin{array}{l} 13,14, \\ 23,34 \end{array}$};
\node (47) [rounded rectangle,draw] at (0,5) {$\begin{array}{l} 12,13, \\ 14,24 \end{array}$};
\node (48) [rounded rectangle,draw] at (1,5) {$\begin{array}{l} 13,14, \\ 24,34 \end{array}$};
\node (49) [rounded rectangle,draw] at (2,5) {$\begin{array}{l} 12,14, \\ 23,24 \end{array}$};
\node (410) [rounded rectangle,draw] at (3,5) {$\begin{array}{l} 12,13, \\ 14,34 \end{array}$};
\node (411) [rounded rectangle,draw] at (4,5) {$\begin{array}{l} 12,13, \\ 14,23 \end{array}$};
\node (412) [rounded rectangle,draw] at (5,5) {$\begin{array}{l} 12,14, \\ 24,34 \end{array}$};

\node[inner sep=7pt,draw,thick,rounded rectangle,rounded rectangle arc length=90,fit=(41) (42) (47) (48)] {};
\node[inner sep=7pt,draw,thick,rounded rectangle,rounded rectangle arc length=90,fit=(43) (44) (49) (410)] {};
\node[inner sep=7pt,draw,thick,rounded rectangle,rounded rectangle arc length=90,fit=(45) (46) (411) (412)] {};
\node[draw,dashed,inner sep=15pt,thick,rounded rectangle,rounded rectangle arc length=120,fit=(43) (44) (49) (410) (45) (46) (411) (412)] {};

\node (31) [rounded rectangle,draw] at (1,4) {$12,13,23$};
\node (32) [rounded rectangle,draw] at (2,4) {$23,24,34$};
\node (33) [rounded rectangle,draw] at (3,4) {$13,14,34$};
\node (34) [rounded rectangle,draw] at (4,4) {$12,14,24$};
\node[inner sep=7pt,draw,thick,rounded rectangle,fit=(31) (32) (33) (34)] {};
\node (35) [rounded rectangle,draw] at (1,3) {$12,23,24$};
\node (36) [rounded rectangle,draw] at (2,3) {$13,23,34$};
\node (37) [rounded rectangle,draw] at (3,3) {$14,24,34$};
\node (38) [rounded rectangle,draw] at (4,3) {$12,13,14$};
\node[inner sep=7pt,draw,thick,rounded rectangle,fit=(35) (36) (37) (38)] {};

\node (21) [rounded rectangle,draw] at (2,2) {$12,34$};
\node (22) [rounded rectangle,draw] at (3,2) {$14,23$};
\node[inner sep=7pt,draw,thick,rounded rectangle,fit=(21) (22)] {};

\node (11) [rounded rectangle,draw] at (1,1) {$12$};
\node (12) [rounded rectangle,draw] at (2,1) {$23$};
\node (13) [rounded rectangle,draw] at (3,1) {$34$};
\node (14) [rounded rectangle,draw] at (4,1) {$14$};
\node[inner sep=7pt,draw,thick,rounded rectangle,fit=(11) (12) (13) (14)] {};

\node (01) [rounded rectangle,draw] at (2.5,0) {$\emptyset$};
\node[inner sep=7pt,draw,thick,rounded rectangle,fit=(01)] {};

\end{tikzpicture}}
\end{center}
\caption{Orbits of the dihedral group action on $\curly{H}$-prime ideals of $\OG q24$.  The $\curly{H}$-prime ideals are labelled by their generating quantum minor labels, groupings indicated by solid lines are orbits under the cycle $c=(1\, 2\, 3\, 4)$ and the grouping indicated by the dashed line is the dihedral orbit that is formed from the two contained cycle orbits.}\label{fig:OqGr24-Hprimes}
\end{figure}

\begin{figure}[ht]
\begin{center} 
\scalebox{0.75}{\input{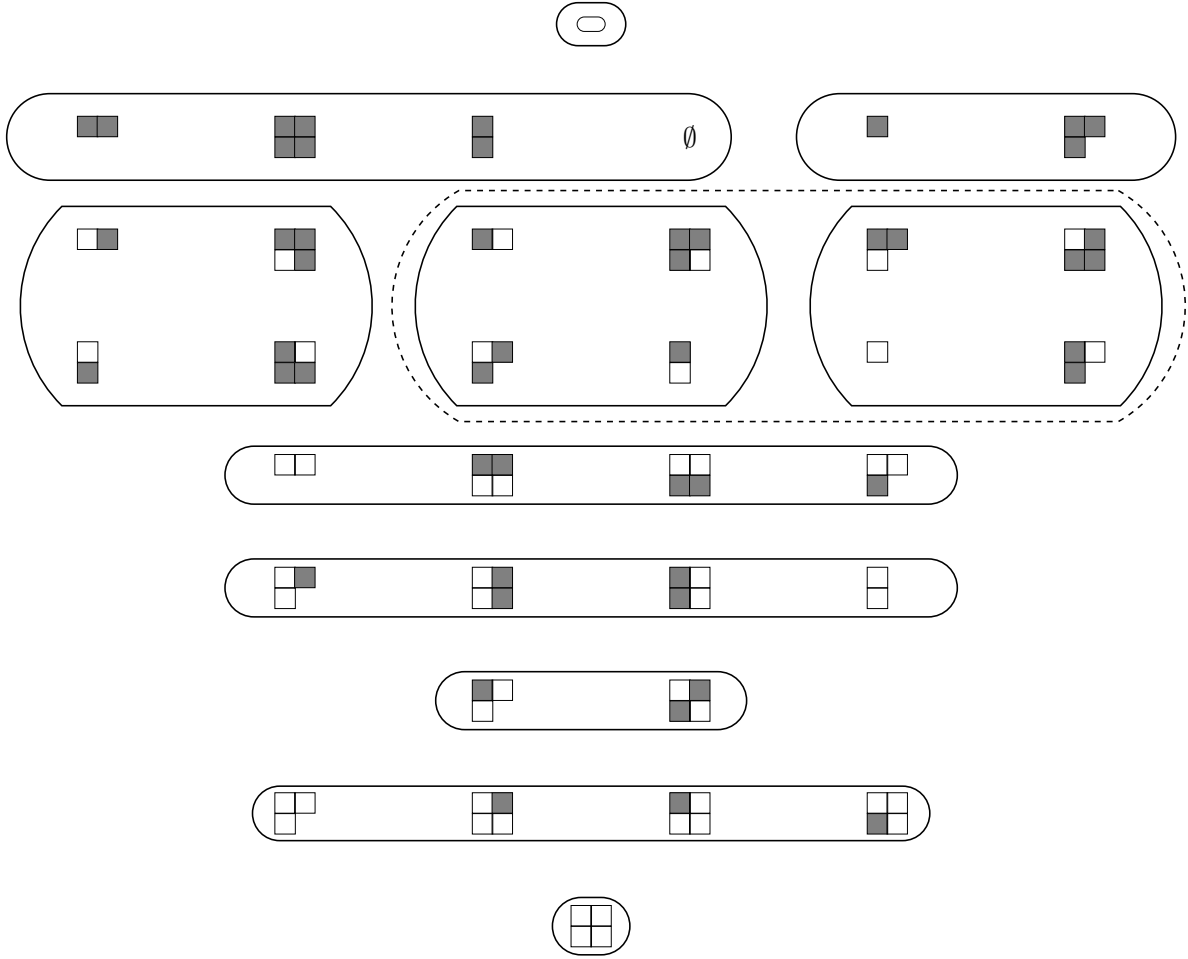}}
\end{center}
\caption{Orbits of the dihedral group action on $\curly{H}$-prime ideals of $\OG q24$.  The $\curly{H}$-prime ideals are labelled by the corresponding Cauchon diagram, groupings indicated by solid lines are orbits under the cycle $c=(1\, 2\, 3\, 4)$ and the grouping indicated by the dashed line is the dihedral orbit that is formed from the two contained cycle orbits.}\label{fig:OqGr24-Hprimes-Cauchon}
\end{figure}

We conclude with an example, that of $\OG q24$.  In this case, it was shown by Russell (\cite{RussellThesis}; see also \cite[Section~2]{LLR}) that every $\curly{H}$-prime ideal is generated by the quantum minors it contains.  (It is still conjectural that this holds for every $\OG qmn$.)  There are 34 $\curly{H}$-prime ideals in $\OG q24$ and as illustrated in Figures~\ref{fig:OqGr24-Hprimes} and \ref{fig:OqGr24-Hprimes-Cauchon} (based on Figure~3 of \cite{LLR}), there are 11 $\lgen c \rgen$-orbits and 10 $I_{2}(n)$-orbits.  We believe it would be valuable to also try to understand this action as an action directly on the set of Cauchon diagrams; some combinatorial features are easily identified but we do not have a complete description of such an action.

\section{A parameterized family of quantum dihedral groupoids}\label{s:groupoid.any.integer}

In \cite{LaunoisLenagan} Launois and Lenagan establish the existence of the isomorphism $\theta$ of Theorem \ref{thm:LL.theta} by considering a certain localisation of $\OG qmn$. In particular, they localise at the consecutive minor with indexing set $M_1\defeq \{ 1,2,\ldots,m\}$ and consider an isomorphism of $\OG qmn [[M_1]^{-1}]$ with $\K[x_{ij}][y^{\pm};\sigma]$, the latter being a skew-Laurent extension of the quantum matrix ring $\K[x_{ij}]$, with $\K[x_{ij}]$ isomorphic to $\OM qm{n-m}$.

However, just as in \cite{LaunoisLenagan}, such a localisation can be carried out at any consecutive minor $[M_\alpha]$ with indexing set $M_{\alpha}=\{ \til\alpha,\til{\alpha+1},\ldots,\til{\alpha+m-1} \}$ and there is an accompanying isomorphism $\OG qmn[[M_\alpha]^{-1}] \isomorphic A_\alpha$ where $A_\alpha$ is again a skew-Laurent extension of a quantum matrix ring isomorphic to $\OM qm{n-m}$.

In what follows, we establish that after twisting $A_\alpha$ by the cocycle $\Gamma$ there is an isomorphism $T(A_\alpha)\isomorphic A_{\alpha+1}$.  We show that the image under this isomorphism of the natural embedding of $T(\OG qmn)$ in $T(A_\alpha)$ is exactly the natural embedding of $\OG qmn$ in $A_{\alpha+1}$. Thus, by restriction, we obtain another isomorphism between $T(\OG qmn)$ and $\OG qmn$, and moreover, we show that on minors it has the effect of cycling indices by $c$, up to a power of $q$. This is completely analogous to the methods of \cite{LaunoisLenagan}, only the map obtained by this procedure acts with different powers of $q$ for each different class of $\alpha$ modulo $n$.

\subsection{Localisation at consecutive minors and dehomogenisation}\label{ss:localisation.consecutive.minors}

For any $\alpha \in \Z$, set $\M\alpha \defeq \{ \til\alpha, \til{\alpha+1}, \ldots , \til{\alpha+m-1} \}$; the associated minor $[M_{\alpha}]$ is called a consec\-utive minor. Then by \cite[Lemma~3.7]{KrobLeclerc}, $[M_{\alpha}]$ quasi-commutes with every element of $\curly{P}_{q}$ and hence is normal in $\OG qmn$.  It follows that one may form the localisation $\locOG \alpha$ in such a way that $\OG qmn$ is embedded as a subring. Define the integers $\dialpha , \djalpha \in \{1,\ldots, n\}$ by
	\BEQ	
	\dialpha \defeq \til{\alpha+m-i}	\;\;\; \mbox{and}	\;\;\;	\djalpha \defeq \til{j+\alpha+m-1}.	\nn
	\EEQ
One useful observation is that for every $\alpha$ and any choice of $i$ and $j$ subject to the constraints $1\leq i \leq m$ and $1\leq j \leq n-m$, we have $\dialpha\in\M\alpha$ and $\djalpha\notin\M\alpha$. In particular, the two never coincide. Now, in the algebra $\locOG \alpha$, set
	\BEQ	
	\x \alpha i j \defeq [\M\alpha \union \{ \djalpha \} \setminus \{ \dialpha \}] \inv{[\M\alpha]}.	\nn
	\EEQ
By \cite[Theorem~4.1]{LaunoisLenagan}, for any $\alpha$ the subalgebra $\K[\x \alpha ij] \subset \locOG \alpha$ is isomorphic to the quantum matrix algebra $\OM{q}{m}{n-m}$ via the map $\x\alpha ij\mapsto X_{ij}$.  Moreover there exists a \emph{dehomogenisation isomorphism} $\phi_{\alpha}$ between the localisation $\locOG{\alpha}$ and a skew-Laurent extension $\dehom{\alpha}$, where $\sigma_\alpha$ is the automorphism of $\K[\x \alpha i j]$ defined by the equalities $\sigma_\alpha(\x \alpha i j) [\M\alpha] = [\M\alpha] \x \alpha ij$ for all $i,j$, with $y_\alpha = \phi_\alpha([\M\alpha])$. We will denote the inverse isomorphism by $\rho_{\alpha}$ so that we have
	\BEQ
	\locOG\alpha \mathop{\rightleftarrows}_{\rho_{\alpha}}^{\phi_{\alpha}} \dehom\alpha. 	\nn
	\EEQ
Throughout the sequel, we will write $\A\alpha \defeq \dehom \alpha$ and we note that $\A\alpha\isomorphic\A\beta$ whenever $\alpha\equiv\beta \mod n$ since in this case $[M_\alpha] = [M_\beta]$. We shall need to understand the structure of the algebras $\A\alpha$ and since the normal element $[\M\alpha]$ quasi-commutes with all of $\curly{P}_{q}$ and hence with each $\x\alpha ij$, we are left only to determine relations of the type $y_\alpha \x \alpha ij = q^\ell \x \alpha ij y_\alpha$ for some $\ell\in\Z$. This is accomplished with the next two results.

\BLEM	\label{lem:A.alpha.structure.case1}
Let $1\leq i \leq m$ and $1\leq j \leq n-m$. Suppose that $\alpha\in\Z$ such that $1 \leq \aalpha \leq n-m$. Then one has
	\begin{align*}
	\sigma_\alpha(\x \alpha ij) &= \BCASE	\inv q\x\alpha ij	&\caseif	j > n-m-\aalpha +1	\\
									q \x \alpha ij 	&\caseif	j \leq n-m-\aalpha +1		\ECASE,	
\intertext{and consequently}
	y_\alpha \x\alpha ij &= \BCASE	\inv q\x\alpha ij y_\alpha	&\caseif	j > n-m-\aalpha +1	\\
							q \x \alpha ij y_\alpha	&\caseif	j \leq n-m-\aalpha +1		\ECASE.	
	\end{align*}
\ELEM

\BPF
Set $N_\alpha \defeq \M\alpha \setminus \{ \dialpha \}$. Then $\x\alpha ij [\M\alpha] = [N_\alpha \union \{ \djalpha \}]$ while $[\M\alpha] = [N_\alpha \union \{ \dialpha \}]$. Note that for integers $1 \leq r < s \leq n$, one has $[r][s]=q[s][r]$ in $\OG q1n$. Using Proposition~\ref{prop:QuantumMuirsLaw}, it follows that $[\M\alpha](\x \alpha ij [\M\alpha]) = q(\x \alpha ij [\M\alpha])[\M\alpha]$ whenever $\dialpha<\djalpha$, whence multiplying by $\inv{[\M\alpha]}$ on the right gives $\sigma_\alpha(\x\alpha ij)=q\x\alpha ij$. Similarly, when $\dialpha>\djalpha$, one obtains that $\sigma_\alpha(\x\alpha ij)=\inv q \x\alpha ij$. To complete the proof, it is a routine calculation to show that the inequality $\dialpha<\djalpha$ holds exactly when $j\leq n-m-\aalpha +1$.
\EPF

\BLEM	\label{lem:A.alpha.structure.case2}
Let $1\leq i \leq m$ and $1\leq j \leq n-m$. Suppose that $\alpha\in\Z$ such that $n-m+1 \leq \aalpha \leq n$. Then one has
	\begin{align*}
	\sigma_\alpha(\x \alpha ij) &= {\BCASE	\inv q\x\alpha ij	&\caseif	i \geq \aalpha - (n-m)	\\
									q \x \alpha ij 	&\caseif	i < \aalpha - (n-m)		\ECASE},		
\intertext{and consequently}
	y_\alpha \x \alpha ij 		&= {\BCASE	\inv q\x\alpha ij y_\alpha	&\caseif	i \geq \aalpha - (n-m)	\\
									q \x \alpha ij y_\alpha 	&\caseif	i < \aalpha - (n-m)		\ECASE}.	
	\end{align*}
\ELEM

\BPF
The proof is similar to that of Lemma \ref{lem:A.alpha.structure.case1} and is omitted.
\EPF

\subsection{Twisting $\A\alpha$}\label{ss:twist.A.alpha}

We now wish to consider the effect of twisting the algebras $\A\alpha$ defined in the previous section. As in Lemmas \ref{lem:A.alpha.structure.case1} and \ref{lem:A.alpha.structure.case2}, the calculations fall into two cases: (i) when $1 \leq \aalpha \leq n-m$ and (ii) when $n-m+1 \leq \aalpha \leq n$. Since $[M_{\alpha}]$ is a homogeneous element of $\OG qmn$, the natural grading of $\OG qmn$ defined by equation (\ref{eqn:natural.grading}) extends to one on $\locOG{\alpha}$ and so the algebra $\A\alpha$ also inherits a $\Z^n$-grading via the isomorphism $\phi_{\alpha}$. In particular, the twist $T(\A\alpha)$ is well defined.

Recall that we are writing $\Gamma(a,b)$ in place of $\Gamma(\content(a),\content(b))$ for elements $a,b\in\A\alpha$ and note that for any $\alpha\in\Z$,
	\BEQ
	\content(\x\alpha ij) = \ep(\djalpha) - \ep(\dialpha)	\hgap\text{and}\hgap	\content(y_\alpha) = \sum_{\nu\in\M\alpha} \ep(\nu).	\nn
	\EEQ
Given this observation, we begin the calculations of the twisted product on pairs from the set $\{ \tx \alpha ij, \hat{y}_\alpha\}\subset T(\A\alpha)$. First we deal with the situation in which $1\leq \aalpha \leq n-m$ for $\alpha\in\Z$.\\

\BLEM	\label{lem:num.theory.case1}
Let $\alpha\in\Z$ such that $1\leq \aalpha \leq n-m$, $1\leq i \leq m$ and $1\leq j \leq n-m$.  Then
	\begin{enumerate}[$(\rm a)$]
	\item $n\notin \M\alpha$,
	\item $\dialpha\leq n-1$; in particular, there is no $i\in\{1,\ldots,m\}$ such that $n=\dialpha$, and
	\item $\djalpha = n$ if and only if $j=n-m-\aalpha+1$.
	\end{enumerate}
\ELEM

\BPF
As an example of the type of calculations involved, we prove (c). 

First note that $j+\aalpha+m-1 \leq 2n-m-1$, so $\djalpha = j+\aalpha + m -1 - \ell n$ where $\ell\in \{0,1\}$. If $\djalpha=n$, then one obtains $j=(1+\ell)n-m-\aalpha+1$. If $\ell=1$, then $j\geq n+1$, a contradiction. So we must have $\ell=0$ and the result follows. Parts (a) and (b) also follow immediately from the bounds on $i$, $j$ and $\aalpha$.
\EPF

As a consequence, we obtain the next result.

\BLEM	\label{lem:c.values.case1}
Let $\alpha\in\Z$ such that $1 \leq \aalpha \leq n-m$ and let $a,b\in\{\x\alpha ij,y_\alpha\}$. Then $\Gamma(a,b)=1$, except for the following two special cases:
	\BEQAR
	\Gamma(\x{\alpha}{i_1,}{n-m-\aalpha+1},\x{\alpha}{i_2,}{n-m-\aalpha+1}) &=& \inv p	\hgap\text{and}	\nn	\\
	\Gamma(\x{\alpha}{i,}{n-m-\aalpha+1}, y_\alpha) &=& q^2,						\nn
	\EEQAR
where $i$, $i_1$ and $i_2$ are elements of the set $\{1,\ldots,m\}$.
\ELEM

\BPF
Note that $\Gamma(a,b)$ can only take a value other than $1$ when $\ep(n)$ appears with a nonzero coefficient in $\content(a)$. Lemma \ref{lem:num.theory.case1} implies that this can only be the case when $a=\x{\alpha}{i,}{n-m-\aalpha+1}$ for some $i\in\{1,\ldots,m\}$. Then, $\Gamma(\x{\alpha}{i,}{n-m-\aalpha+1}, y_\alpha)=p^{(1)(m)}=q^2$ since $\ep(n)$ appears with coefficient $1$ in $\content(\x{\alpha}{i,}{n-m-\aalpha+1})$ but does not appear in $\content(y_\alpha)$. The other calculations are similar and therefore omitted.
\EPF

The next two lemmas provide the analogous results when $n-m+1\leq\aalpha\leq n$.

\BLEM	\label{lem:num.theory.case2}
Let $\alpha\in\Z$ such that $n-m+1\leq \aalpha \leq n$, $1\leq i \leq m$ and $1\leq j \leq n-m$.  Then
	\begin{enumerate}[$(\rm a)$]
	\item $n\in\M\alpha$,
	\item $\dialpha =n$ if and only if $i=\aalpha-(n-m)$, and
	\item there is no $j\in\{1,\ldots,n-m\}$ such that $n=\djalpha$.
	\end{enumerate}
\ELEM

\BPF
We prove (c). Since $n+1 \leq j+\aalpha+m-1\leq 2n-1$, we have that $\djalpha = j+\aalpha + m -1 - n$. Thus, if $n=\djalpha$, one obtains that $j=2n-m-\aalpha+1\geq n-m+1$, a contradiction. As in Lemma \ref{lem:num.theory.case1}, parts (a) and (b) follow similarly from the bounds on $i$, $j$ and $\aalpha$.
\EPF

\BLEM	\label{lem:c.values.case2}
Let $\alpha\in\Z$ such that $n-m+1 \leq \aalpha \leq n$ and let $a,b\in\{\x\alpha ij,y_\alpha\}$. Then $\Gamma(a,b)=1$, except for the following three special cases:
	\BEQAR
	\Gamma(\x{\alpha}{\aalpha-(n-m),}{j_1},\x{\alpha}{\aalpha-(n-m),}{j_2}) 	&=& \inv p,				\nn	\\
	\Gamma(\x{\alpha}{\aalpha-(n-m),}{j},y_\alpha) 						&=& pq^{-2},	\hgap\text{and}	\nn	\\
	\Gamma(y_\alpha, \x{\alpha}{\aalpha-(n-m),}{j}) 						&=& p,					\nn
	\EEQAR
where $j$, $j_1$ and $j_2$ are elements of the set $\{1,\ldots,n-m\}$.
\ELEM

\BPF
The proof is similar to that of Lemma \ref{lem:c.values.case1}, only using Lemma \ref{lem:num.theory.case2} to give the conditions under which $\ep(n)$ appears with a nonzero coefficient in $\content(a)$.
\EPF

\subsection{The isomorphism of $T(\A\alpha)$ with $\A{\alpha+1}$}\label{ss:isomforAalphas}

We are now in a position to determine the relations in the algebra $T(\A\alpha)$. Since $\A\alpha$ is generated by the homogeneous elements $\{\x\alpha ij,y_\alpha\}$ then $T(\A\alpha)$ is generated by the homogeneous elements $\{\Tx\alpha ij,T(y_{\alpha})\}$.

\BLEM	\label{lem:structure.K[x']}
For every $\alpha\in\Z$, $(\tx\alpha ij)$ is a generic $q$-quantum matrix; that is, the algebra $\K[\tx\alpha ij]$ is isomorphic to $\OM{q}{m}{n-m}$.
\ELEM

\BPF
The verifications are trivial in the cases for which $\Gamma(-,-)=1$. Otherwise, there are two situations to consider: (i) when $1\leq \aalpha \leq n-m$ and (ii) when $n-m+1\leq \aalpha \leq n$. First consider (i). Suppose that $i_1<i_2$. By Lemma \ref{lem:c.values.case1} the only nontrivial computation that remains is
	\begin{align*}
	\Tx\alpha{i_1,}{n-m-\aalpha+1}\Tx\alpha{i_2,}{n-m-\aalpha+1}
	&=	\Gamma(\x{\alpha}{i_1,}{n-m-\aalpha+1},\x{\alpha}{i_2,}{n-m-\aalpha+1})
				T(\x\alpha{i_1,}{n-m-\aalpha+1}\x\alpha{i_2,}{n-m-\aalpha+1})	\\
	&=	\inv p T(q\,\x\alpha{i_2,}{n-m-\aalpha+1}\x\alpha{i_1,}{n-m-\aalpha+1})	
\intertext{On the other hand,}
	\Tx\alpha{i_2,}{n-m-\aalpha+1}\Tx\alpha{i_1,}{n-m-\aalpha+1}
	&=	\Gamma(\x{\alpha}{i_2,}{n-m-\aalpha+1},\x{\alpha}{i_1,}{n-m-\aalpha+1})
				T(\x\alpha{i_2,}{n-m-\aalpha+1}\x\alpha{i_1,}{n-m-\aalpha+1})	\\
	&=	\inv p T(\x\alpha{i_2,}{n-m-\aalpha+1}\x\alpha{i_1,}{n-m-\aalpha+1}).
\end{align*}	 
Therefore
\begin{align*}
 \Tx\alpha{i_1,}{n-m-\aalpha+1}\Tx\alpha{i_2,}{n-m-\aalpha+1} & = q \Tx\alpha{i_2,}{n-m-\aalpha+1}\Tx\alpha{i_1,}{n-m-\aalpha+1},
\intertext{or equivalently} \tx\alpha{i_1,}{n-m-\aalpha+1}\tx\alpha{i_2,}{n-m-\aalpha+1} & = q \tx\alpha{i_2,}{n-m-\aalpha+1}\tx\alpha{i_1,}{n-m-\aalpha+1} 
\end{align*}
in the ``hat'' notation, as desired. In case (ii), the only the interesting case is the row relation corresponding to $i=\aalpha-(n-m)$. The calculation is similar to the one above and is omitted. 
\EPF

To have a complete set of relations for $T(\A\alpha)$, it only remains to compute commutation relations between $\hat{y}_\alpha$ and the elements $\tx\alpha ij$. As usual, it is helpful to consider the calculations in two cases.

\BLEM	\label{lem:T.A.alpha.structure.case1}
Let $\alpha\in\Z$ such that $1\leq \aalpha \leq n-m$.
	\BEQ
	\hat{y}_\alpha \tx\alpha ij = \BCASE	\inv q\tx\alpha ij \hat{y}_\alpha	&\caseif	j > n-m-\aalpha		\\
								q \tx\alpha ij \hat{y}_\alpha	&\caseif	j \leq n-m-\aalpha	\ECASE.	\nn
	\EEQ
\ELEM

\BPF
We use Lemma \ref{lem:A.alpha.structure.case1} for commutation relations in $\A\alpha$. According to Lemma \ref{lem:c.values.case1}, the only nontrivial case (i.e.\ when $\Gamma(-,-)\neq 1$) occurs when $j=n-m-\aalpha+1$. We have for any $i\in\{1,\ldots,m\}$ that
	\begin{align*}
	T(y_\alpha) \Tx\alpha{i,}{n-m-\aalpha+1}
	&=	\Gamma(y_\alpha, \x\alpha{i,}{n-m-\aalpha+1})T(y_\alpha \x\alpha{i,}{n-m-\aalpha+1})		\\
	&=	T(y_\alpha \x\alpha{i,}{n-m-\aalpha+1}) = qT(\x\alpha{i,}{n-m-\aalpha+1} y_\alpha).	
\intertext{Alternatively,}
\Tx\alpha{i,}{n-m-\aalpha+1} T({y}_\alpha)
	&=	\Gamma(\x\alpha{i,}{n-m-\aalpha+1} , y_\alpha) T(\x\alpha{i,}{n-m-\aalpha+1} y_\alpha) 		\\
	&=	q^2 T(\x\alpha{i,}{n-m-\aalpha+1} y_\alpha).								
	\end{align*}
Thus, reverting to ``hat'' notation, $\hat{y}_\alpha \tx\alpha{i,}{n-m-\aalpha+1} = \inv q \, \tx\alpha{i,}{n-m-\aalpha+1} \hat{y}_\alpha$ and the result follows.
\EPF

\BLEM	\label{lem:T.A.alpha.structure.case2}
Let $\alpha\in\Z$ such that $n-m+1\leq \aalpha \leq n$. Then
	\BEQ
	\hat{y}_\alpha \tx\alpha ij = \BCASE	\inv q\tx\alpha ij \hat{y}_\alpha	&\caseif	i \geq \aalpha +1 - (n-m)	\\
								q \tx\alpha ij \hat{y}_\alpha	&\caseif	i < \aalpha +1 - (n-m) 	\ECASE.	\nn
	\EEQ
\ELEM

\BPF
By Lemma \ref{lem:c.values.case2}, the only nontrivial computation occurs when $i=\aalpha-(n-m)$. The rest follows directly from Lemma \ref{lem:A.alpha.structure.case2}. Consider that
	\begin{align*}
	T(y_\alpha) \Tx\alpha{\aalpha-(n-m),}{j}
	&=	\Gamma(y_\alpha, \x\alpha{\aalpha-(n-m),}{j}) T(y_\alpha \x\alpha{\aalpha-(n-m),}{j})		\\
	&=	p T(\inv q \x\alpha{\aalpha-(n-m),}{j} y_\alpha).							
\intertext{On the other hand}
	\Tx\alpha{\aalpha-(n-m),}{j}T({y}_\alpha)
	& = \Gamma(\x\alpha{\aalpha-(n-m),}{j}, y_\alpha) T(\x\alpha{\aalpha-(n-m),}{j} y_\alpha) 		\\
	& = pq^{-2} T(\x\alpha{\aalpha-(n-m),}{j} y_\alpha)	.							
	\end{align*}
Hence $\hat{y}_\alpha \tx\alpha{\aalpha-(n-m),}{j} = q\, \tx\alpha{\aalpha-(n-m),}{j} \hat{y}_\alpha$ and the claim is proved.
\EPF

These results lead to the following theorem, which relates the twist of the algebra $\A\alpha$ to the algebra $\A{\alpha+1}$.

\BTHM	\label{thm:T.A.alpha.isomorphic.A.alpha+1}
For every integer $\alpha$, there exists an isomorphism
	\BEQ
	\theta_\alpha \colon T(\A\alpha) \isomap \A{\alpha+1}	\nn
	\EEQ
sending $\tx\alpha ij \mapsto \x{\alpha+1}ij$ and $\hat{y}_\alpha \mapsto y_{\alpha+1}$.
\ETHM

\BPF
We begin by showing that the relations among the respective generating sets $\{\tx\alpha ij, \hat{y}_\alpha\}$ and $\{\x{\alpha+1} ij, y_{\alpha+1}\}$ coincide.  Using Lemmas \ref{lem:A.alpha.structure.case1}, \ref{lem:A.alpha.structure.case2}, \ref{lem:structure.K[x']}, \ref{lem:T.A.alpha.structure.case1} and \ref{lem:T.A.alpha.structure.case2}, the verification is straightforward except for the two boundary cases $\aalpha=n-m$ and $\aalpha=n$. Consider the case $\aalpha=n-m$. Since both $(\tx\alpha ij)$ and $(\x{\alpha+1}ij)$ are generic $q$-quantum matrices we need only to compare the commutation relations of the $x_{ij}$'s with the relevant $y$'s. Lemma \ref{lem:T.A.alpha.structure.case1} implies that when $\aalpha=n-m$, one has $\hat{y}_\alpha \tx\alpha ij = \inv q \tx\alpha ij \hat{y}_\alpha$ for every pair $(i,j)\in\{1,\ldots,m\}\times\{1,\ldots,n-m\}$. On the other hand, Lemma \ref{lem:A.alpha.structure.case2} implies that when $\aalpha=n-m+1$, one has $y_{\alpha+1}\x{\alpha+1}ij = \inv q \x{\alpha+1}ij y_{\alpha+1}$ for any pair $(i,j)$, as desired. The proof for $\aalpha=n$ is similar.

Now we may define an onto homomorphism $\theta_{\alpha}$ as above.  Then since $\K[\tx \alpha ij]$ and $\K[\x {\alpha+1}{i}{j}]$ are isomorphic to $\OM{q}{m}{n-m}$, both have Gelfand--Kirillov dimension $m(n-m)$.  It follows that $\A\alpha$ and $\A{\alpha+1}$ have Gelfand--Kirillov dimension $m(n-m)+1$ and so does $T(\A\alpha)$, since there is a graded isomorphism between $\A\alpha$ and $T(\A\alpha)$.  Furthermore, $\A\alpha$ and $\A{\alpha+1}$ are domains and twisting by a 2-cocycle preserves the property of being an integral domain (\cite[Lemma~5.2]{LaunoisLenagan}) and so $T(\A\alpha)$ is also a domain.  Since any epimorphism between domains of the same Gelfand--Kirillov dimension is an isomorphism (\cite[Proposition~3.15]{KrauseLenagan}), the result follows.
\EPF

\subsection{Twisting $\OG qmn$}\label{ss:isomforOqGr}

One may regard $\OG qmn$ as a subalgebra of $\A\beta$, embedded via the map $\phi_\beta$, following the natural inclusion $\OG qmn \hookrightarrow \locOG\beta$. In the context of Theorem \ref{thm:T.A.alpha.isomorphic.A.alpha+1} there are then two relevant embeddings of $\OG qmn$ for each $\alpha\in\Z$, a twisted version $T(\OG qmn) \subset T(\A\alpha)$ and another copy $\OG qmn \subset \A{\alpha+1}$. Our next goal is to show that these two algebras correspond under the isomorphism $\theta_\alpha$.  To do this, one must trace the image of a generating quantum minor in $\OG qmn$ under the sequence of maps
	\BEQ
	\OG qmn\hookrightarrow\locOG\alpha	\mapname{\phi_\alpha}	\A\alpha 		\mapname{T}
	T(\A\alpha)						\mapname{\theta_\alpha} 	\A{\alpha+1}	\mapname{\rho_{\alpha+1}}	\locOG{\alpha+1}.	\nn
	\EEQ
Recall that as a notational convenience, for $j\in\Z$ we write $\til j$ to mean the representative from $\{1,2,\ldots,n\}$ which is congruent to $j$ modulo $n$. Moreover, if $J\subset\{1,\ldots,n\}$ and $k\in\Z$, we let $J + k = k+J = \{\til{j + k} \st j\in J\}$.

Given a row set $K=\{k_1,\ldots,k_t\}$ with $1\leq k_1< \cdots < k_t \leq m$ and a column set $L=\{l_1,\ldots,l_t\}$ with $1\leq l_1<\ldots<l_t\leq n-m$, we denote the corresponding quantum minor in $\A\alpha$ (respectively $T(\A\alpha)$) by $\qminor KL\alpha$ (resp.\ $\tqminor KL\alpha$). We begin with the following lemma, which is a straightforward generalisation of \cite[Lemma~5.6]{LaunoisLenagan}.

\BLEM	\label{lem:twist.effect.quantum.minors}
Let $K$ and $L$ be as above, so that $\qminor KL\alpha \in \A\alpha$. Then $T(\qminor KL\alpha) = \tqminor KL \alpha$. \qed
\ELEM

We will also need to know the values taken by the cocycle $\Gamma$ on general quantum minors.

\BLEM	\label{lem:c.value.quantum.minors}
Let $K$ and $L$ be as above, so that $\qminor KL\alpha \in \A\alpha$.
	\begin{enumerate}[$(\rm a)$]
	\item If $1\leq \aalpha \leq n-m$, then
		\BEQ
		\Gamma(\qminor KL\alpha , y_\alpha) = 
		\BCASE
			q^2	&\caseif	n-m-\aalpha+1\in L	\\
			1	&	\text{otherwise}
		\ECASE.	\nn
		\EEQ
	\item If $n-m+1\leq \aalpha \leq n$, then
		\BEQ
		\Gamma(\qminor KL\alpha , y_\alpha) = 
		\BCASE
			pq^{-2}	&\caseif	\aalpha-(n-m)\in K	\\
			1		&	\text{otherwise}
		\ECASE.	\nn
		\EEQ
	\end{enumerate}
\ELEM

\BPF
As a result of \cite[Proposition~4.3]{LenaganRussell}, we immediately have
	\BEQ
	\content \left(\qminor KL\alpha \right) = \sum_{\nu=1}^t \left(\ep(\til{l_\nu+\alpha+m-1})-\ep(\til{\alpha+m-k_\nu})\right).	\nn
	\EEQ
Note that in the summation on the right-hand side the first argument is the integer $\djalpha$ for $j=l_\nu$, while the second is $\dialpha$ for $i=k_\nu$. Hence, we may use Lemmas \ref{lem:num.theory.case1} and \ref{lem:num.theory.case2} to determine when $\ep(n)$ appears with a nonzero coefficient in $\content(\qminor KL\alpha)$. For example, when $n-m+1\leq \aalpha \leq n$, Lemma \ref{lem:num.theory.case1}(b) gives that $\Gamma(\qminor KL\alpha , y_\alpha)\neq 1$ only if $\aalpha-(n-m)\in K$. Now Lemma \ref{lem:num.theory.case2}(a) implies that $\ep(n)$ always appears in $\content(y_\alpha)$ with coefficient $1$, so that if $\aalpha-(n-m)\in K$, we obtain
	\BEQ
	\Gamma(\qminor KL\alpha , y_\alpha) = p^{(-1)(m-1)}=pq^{-2},	\nn
	\EEQ
as desired. The case $1\leq\aalpha\leq n-m$ is similar.
\EPF

The following proposition is then the key to obtaining the result we seek.

\BPROP	\label{lem:image.generating.quantum.minor}
Let $[I] \in \OG qmn \subset \locOG\alpha$ as above. Then
	\BEQ
	\rho_{\alpha+1} \compose \theta_\alpha \compose T \compose \phi_\alpha([I]) = \lambda_{I}^\alpha \multdot [I+1],	\nn
	\EEQ
where
	\BEQ
	\lambda_I^\alpha = \BCASE
	q^{-2}	&\caseif	1\leq \aalpha \leq n-m \; \text{and} \; n\in I			\\
	q^2/p	&\caseif	n-m+1\leq \aalpha \leq n \; \text{and} \; n\notin I	\\
	1		& \text{otherwise}
	\ECASE.	\nn
	\EEQ
\EPROP

\BPF
The proof is essentially that of \cite[Lemma~5.8]{LaunoisLenagan}, but with a few modifications for the more general setting. Fix $\alpha\in\Z$ and write $I_r \defeq I\intersect\M\alpha$ and $I_c \defeq I\setminus I_r$. By \cite[Corollary~4.2]{LaunoisLenagan} and \cite[Proposition~4.3]{LenaganRussell}, we have the formul\ae
\begin{align}
\label{eqn:phi.formula} \phi_\alpha([I])	& = \qminor{(\aalpha+m) - (\M\alpha\setminus I_r)}{I_c - (\aalpha+m-1)}\alpha  y_\alpha	\\ 
\intertext{and}
\label{eqn:rho.formula} \rho_\alpha(\qminor KL\alpha) & = [\M\alpha\setminus((\aalpha+m)-K) \disjunion ((\aalpha+m-1)+L)] \inv{[\M\alpha]}		
\end{align}
where $K$ and $L$ are as above. Moreover, for any $\alpha$, one has the identity of sets
	\BEQAR
	\lefteqn{\qminor{(m+\aalpha)-(\M\alpha\setminus I_r)}{I_c-(\aalpha+m-1)}{\alpha+1} = }		\nn \\
	& & \qminor{(m+\aalpha+1)-(\M{\alpha+1}\setminus(I_r+1))}{(I_c+1)-(\aalpha+m)}{\alpha+1}.	\label{eqn:cycle.sets}
	\EEQAR
We will do the computations in the case $n-m+1\leq\aalpha\leq n$; the case $1\leq\aalpha\leq n-m$ is similar. By equation (\ref{eqn:phi.formula}) and Lemma \ref{lem:twist.effect.quantum.minors} we have
	\BEQAR
	T\compose \phi_\alpha([I])
	&=&	T\left( \qminor{(\aalpha+m) - (\M\alpha\setminus I_r)}{I_c - (\aalpha+m-1)}\alpha  y_\alpha \right)	\nn	\\
	&=&	\inv\mu \tqminor{(\aalpha+m) - (\M\alpha\setminus I_r)}{I_c - (\aalpha+m-1)}\alpha \hat{y}_\alpha,		\nn					
	\EEQAR
where $\mu \defeq \Gamma(\qminor{(\aalpha+m) - (\M\alpha\setminus I_r)}{I_c - (\aalpha+m-1)}\alpha,  y_\alpha)$. Therefore, Lemma \ref{lem:c.value.quantum.minors}(b) implies that $\mu=pq^{-2}$ whenever $\aalpha-(n-m) \in (\aalpha+m) - (\M\alpha\setminus I_r)$, and $\mu=1$ otherwise. That is, $\mu=pq^{-2}$ exactly when
	\BEQ
	\aalpha-(n-m) \in (\aalpha+m) - (\M\alpha\setminus I_r)	\iff	n\in\M\alpha\setminus I_r	\iff	n \notin I,		\nn
	\EEQ
where the last implication follows from Lemma \ref{lem:num.theory.case2}(a). Applying $\theta_\alpha$, we obtain
	\BEQAR
	\theta_\alpha \compose T \compose \phi_\alpha([I])
	&=&	\inv\mu \qminor{(m+\aalpha)-(\M\alpha\setminus I_r)}{I_c-(\aalpha+m-1)}{\alpha+1} y_{\alpha+1}			\nn	\\
	&=&	\inv\mu \qminor{(m+\aalpha+1)-(\M{\alpha+1}\setminus(I_r+1))}{(I_c+1)-(\aalpha+m)}{\alpha+1} y_{\alpha+1},	\nn
	\EEQAR
where the last equality is given by equation (\ref{eqn:cycle.sets}). Finally, we use equation (\ref{eqn:rho.formula}) to see that
	\BEQAR
	\rho_{\alpha+1} \compose \theta_\alpha \compose T \compose \phi_\alpha([I])
	&=&	\inv\mu [\M{\alpha+1}\setminus((m+\aalpha+1)-(\M{\alpha+1}\setminus(I_r+1)) \disjunion		\nn	\\
	&~&	\hgap \hgap \hgap \hgap	((\aalpha+m)+(I_c+1)-(\aalpha+m))] \inv{[\M\alpha]} [\M\alpha]		\nn	\\
	&=&	\inv\mu [(\M{\alpha+1}\setminus\{\M{\alpha+1}\setminus(I_r+1)\}) \disjunion (I_c+1)]			\nn	\\
	&=&	\inv\mu [(I_r+1) \disjunion (I_c+1)] \nn \\
        &=& \inv\mu [I+1].									\nn
	\EEQAR
Substituting the appropriate value of $\mu$ proves the claim.
\EPF

The previous proposition allows us to conclude

\BTHM
For each $\alpha\in\Z$, there exists an isomorphism $T(\OG qmn) \isomap \OG qmn$ which sends $T([I])$ to $\lambda_I^\alpha \multdot [{I+1}]$.
\ETHM

\BPF
Let $\iota_\beta:\OG qmn \hookrightarrow \locOG\beta$ denote the natural embedding. Now, identify the source $T(\OG qmn)$ with the image of $\OG qmn$ under the composition $T\compose \phi_\alpha\compose \iota_\alpha$. Similarly, identify the target $\OG qmn$ with its image under $\iota_{\alpha+1}$. Then the isomorphism described above is given by the composition of isomorphisms $\rho_{\alpha+1}\compose\theta_\alpha$ as per the previous proposition.
\EPF

\BRMK
Observe that when $\alpha=1$, we obtain exactly Theorem \ref{thm:LL.theta}; that is, the result of Launois and Lenagan \cite[Theorem~5.9]{LaunoisLenagan}. So, their result should be viewed as a member of a parameterized family of similar results, one for each value of $\alpha$ modulo $n$. Moreover, note that the construction of the dihedral groupoid in Section \ref{ss:groupoid.1} goes through for any $\alpha$, with the same maps $\Omega_{\ell}$ as ``quantum reflections''. So in fact there exists a parameterized family of groupoids, each of which has a graphical representation resembling that on page~\pageref{groupoid}. The arrows, just as in Sections \ref{ss:groupoid.1} and \ref{ss:dihedral}, act like the dihedral group (up to a power of $q$), only the power of $q$ is different for each class $\tilde\alpha$ modulo $n$.

In fact this is indicative of something deeper. Whenever one can prove a result about the quantum Grassmannian by way of localisation at a consecutive minor and the dehomogenisation isomorphism, the result should exist in a parameterized family of similar results, one for each consecutive minor.  Another example of this phenomenon may be seen in \cite{LenaganRussell}, where a parameterized family of partial orders lead to different structures of a graded quantum algebra with a straightening law on $\OG qmn$.
\ERMK

\section{A dihedral action on the totally nonnegative Grassmannian}\label{s:tnn}

We conclude with the observation that the dihedral action we have studied above also exists on the totally nonnegative and totally positive Grassmannians.  In this section, we specialise to considering $\K = \mathbb{R}$, so that all matrices will have real entries.  We denote by $M^{+}(m,n)$ the space of real $m\cross n$ matrices having all $m\cross m$ minors nonnegative and let $\mathrm{GL}^{+}(m)$ be the group of all real $m\cross m$ matrices with positive determinant.  The group $\mathrm{GL}^{+}(m)$ acts naturally on $M^{+}(m,n)$ by left multiplication.

\BDEFN
Define the totally nonnegative Grassmannian $G^{\mathrm{tnn}}(m,n)$ to be the quotient space $\mathrm{GL}^{+}(m)\setminus M^{+}(m,n)$.  The totally positive Grassmannian $G^{\mathrm{tp}}(m,n)$ is the subspace of $G^{\mathrm{tnn}}(m,n)$ of matrices having all $m\cross m$ minors strictly positive.
\EDEFN

Note that $G^{\mathrm{tnn}}(m,n)$ is a closed subset of $G(m,n)$ and $G^{\mathrm{tp}}(m,n)$ is an open subset, both of the same dimension, $k(n-k)$, as $G(m,n)$.  By examining which minors of an element of $G^{\mathrm{tnn}}(m,n)$ are zero and which are positive, one may attribute each element to a cell, giving a cellular decomposition of $G^{\mathrm{tnn}}(m,n)$ into totally nonnegative cells.  The top-dimensional cell of $G^{\mathrm{tnn}}(m,n)$ is $G^{\mathrm{tp}}(m,n)$, for example.

Postnikov (\cite{Postnikov}) has shown that the totally nonnegative cells are in bijection with the so-called Le-diagrams, which have also appeared under the name Cauchon diagrams, and as noted above, by the results of Launois, Lenagan and Rigal in \cite{LLR}, these are in bijection with the $\curly{H}$-prime ideals of $\OG qmn$ except for the augmentation ideal.  This has been studied in more detail in work of Goodearl, Launois and Lenagan (\cite{GLL-TNN-Cells},\cite{GLL-TorusInvPrimes}) and is described in a recent survey by Launois and Lenagan (\cite{LL-TNN-QM-Poisson}).  Thus the totally nonnegative Grassmannian has close links with the quantum Grassmannian.

Postnikov has also observed that one may define an action of a cyclic group of order $n$ on $G^{\mathrm{tnn}}(m,n)$ (\cite[Remark~3.3]{Postnikov}).  We now recall this and show how it may be extended to a dihedral group action.

As above, let $c=(1\, 2\, 3\, \cdots \, n)$ and let $w_{0}$ be the longest element in $S_{n}$, so that $w_{0}$ is the bijection sending $i$ to $n-i+1$ for all $i\in \{1,\ldots,n\}$.  Define $R_{m}\in M(m,m)$ by
\[ (R_{m})_{ij}=\begin{cases} (-1)^{[m/2]} & \text{if}\ i=j=1 \\ 1 & \text{if}\ i=j\neq 1 \\ 0 & \text{otherwise} \end{cases} \]
where $[i]$ denotes the integer part of $i$.

 For a matrix $A\in G^{\mathrm{tnn}}(m,n)$ given as $A=(v_{1}\, v_{2}\, \cdots \, v_{n})$ with columns $v_{i}\in \mathbb{R}^{m}$, define the following actions of $c$ and $w_{0}$ on $A$:
\begin{align*}
c\cdot A & \defeq \left( (-1)^{m-1}v_{n}\, v_{1}\, v_{2} \, \cdots \, v_{n-2} \, v_{n-1} \right),\\
w_{0} \cdot A & \defeq R_{m} \left(v_{n} \, v_{n-1} \, v_{n-2} \, \cdots \, v_{2} \, v_{1}\right).
\end{align*}
That is, $c$ acts by cycling the columns of $A$ and multiplying the first column of the resulting matrix by $(-1)^{m-1}$ and $w_{0}$ acts by permuting the columns of $A$ in the natural way and multiplying the first row of the resulting matrix by $(-1)^{[m/2]}$.  (Note that $(-1)^{m-1}$ and $(-1)^{[m/2]}$ are the signs of the permutations $c$ and $w_{0}$, respectively.)

It is straightforward to see that if $\Delta_{I}(A)$ is an $m\cross m$ minor of $A$ defined by the choice of column subset $I=\{ i_{1} < i_{2} < \cdots < i_{m} \}$, then $\Delta_{I}(A)=\Delta_{c(I)}(c \cdot A)=\Delta_{w_{0}(I)}(w_{0} \cdot A)$, where $c$ and $w_{0}$ act on $I$ as permutations in the natural way.  These hold because the signs introduced in the definitions of the action precisely correct for the signs arising when permuting the columns under the determinants.  So both $c$ and $w_{0}$ preserve total nonnegativity and total positivity.  (We note that acting $n$ times by $c$ on $A$ yields $(-1)^{m-1}A$, which is not equal to $A$ if $m$ is even.  However since $(-1)^{m-1}A$ has the same row space as $A$, it determines the same point in $G(m,n)$.)

Then we see that the dihedral relations hold as follows:
\begin{align*} (w_{0} c w_{0})\cdot A & =(w_{0} c) \cdot R_{m}(v_{n} \, v_{n-1} \, v_{n-2} \, \cdots \, v_{2} \, v_{1}) \\
 & = w_{0} \cdot R_{m}((-1)^{m-1}v_{1} \, v_{n} \, v_{n-1} \, \cdots \, v_{3} \, v_{2}) \\
& = R_{m}^{2}(v_{2} \, v_{3} \, \cdots \, v_{n-1} \, v_{n} \, (-1)^{m-1}v_{1}) \\
& = c^{-1} \cdot A.
\end{align*}
Hence we have an action of the subgroup $\lgen c,w_{0} \rgen \subseteq S_{n}$ on $G^{\mathrm{tnn}}(m,n)$ and $G^{\mathrm{tp}}(m,n)$ and thus the desired dihedral action on the totally nonnegative and totally positive Grassmannians.  As noted in the introduction, it would be interesting to know the implications of this action for the geometric structure of $G^{\mathrm{tnn}}(m,n)$ and $G^{\mathrm{tp}}(m,n)$, particularly their stratifications.
 
\small

\bibliographystyle{amsplain}
\bibliography{QDA-references}\label{references}

\providecommand{\bysame}{\leavevmode\hbox to3em{\hrulefill}\thinspace}
\providecommand{\MR}{\relax\ifhmode\unskip\space\fi MR }
\providecommand{\MRhref}[2]{%
  \href{http://www.ams.org/mathscinet-getitem?mr=#1}{#2}
}
\providecommand{\href}[2]{#2}
\begin{thebibliography}{10}

\bibitem{AsScSh-clusterautom}
Ibrahim Assem, Ralf Schiffler, and Vasilisa Shramchenko, \emph{Cluster
  automorphisms}, Proc. London Math. Soc., to appear.

\bibitem{Brown-Goodearl}
Ken~A. Brown and Ken~R. Goodearl, \emph{Lectures on algebraic quantum groups},
  Advanced Courses in Mathematics. CRM Barcelona, Birkh\"auser Verlag, Basel,
  2002.

\bibitem{FZ-CA2}
Sergey Fomin and Andrei Zelevinsky, \emph{Cluster algebras. {II}. {F}inite type
  classification}, Invent. Math. \textbf{154} (2003), no.~1, 63--121.

\bibitem{GLL-TorusInvPrimes}
Ken~R. Goodearl, St\'{e}phane Launois, and Tom~H. Lenagan,
  \emph{Torus-invariant prime ideals in quantum matrices, totally nonnegative
  cells and symplectic leaves}, Math. Z. \textbf{269} (2011), no.~1, 29--45.

\bibitem{GLL-TNN-Cells}
\bysame, \emph{Totally nonnegative cells and matrix {P}oisson varieties}, Adv.
  Math. \textbf{226} (2011), 779--826.

\bibitem{Gr2nSchubertQCA}
Jan~E. Grabowski and St\'{e}phane Launois, \emph{Quantum cluster algebra
  structures on quantum {G}rassmannians and their quantum {S}chubert cells: the
  finite-type cases}, Int. Math. Res. Not. (2011), no.~10, 2230--2262.

\bibitem{KLR}
Ann~C. Kelly, Tom~H. Lenagan, and Laurent Rigal, \emph{Ring theoretic
  properties of quantum {G}rassmannians}, J. Algebra Appl. \textbf{3} (2004),
  no.~1, 9--30.

\bibitem{KnutsonLamSpeyer}
Allen Knutson, Thomas Lam, and David~E. Speyer, \emph{Positroid varieties:
  juggling and geometry}, preprint, arXiv:1111.3660.

\bibitem{KrauseLenagan}
G{\"u}nter~R. Krause and Tom~H. Lenagan, \emph{Growth of algebras and
  {G}elfand-{K}irillov dimension}, revised ed., Graduate Studies in
  Mathematics, vol.~22, American Mathematical Society, Providence, RI, 2000.

\bibitem{KrobLeclerc}
Daniel Krob and Bernard Leclerc, \emph{Minor identities for quasi-determinants
  and quantum determinants}, Comm. Math. Phys. \textbf{169} (1995), no.~1,
  1--23.

\bibitem{LL-TNN-QM-Poisson}
St\'{e}phane Launois and Tom~H. Lenagan, \emph{From totally nonnegative
  matrices to quantum matrices and back, via {P}oisson geometry}, to appear in
  the Proceedings of the Belfast Workshop on Algebra, Combinatorics and
  Dynamics 2009.

\bibitem{LaunoisLenagan}
\bysame, \emph{Twisting the quantum grassmannian}, Proc. Amer. Math. Soc.
  \textbf{139} (2011), 99--110.

\bibitem{LLR}
St\'{e}phane Launois, Tom~H. Lenagan, and Laurent Rigal, \emph{Prime ideals in
  the quantum {G}rassmannian}, Selecta Math. (N.S.) \textbf{13} (2008), no.~4,
  697--725.

\bibitem{LeclercZelevinsky}
Bernard Leclerc and Andrei Zelevinsky, \emph{Quasicommuting families of quantum
  {P}l\"ucker coordinates}, Kirillov's seminar on representation theory, Amer.
  Math. Soc. Transl. Ser. 2, vol. 181, Amer. Math. Soc., Providence, RI, 1998,
  pp.~85--108.

\bibitem{LR-qGrAlgs}
Tom~H. Lenagan and Laurent Rigal, \emph{Quantum graded algebras with a
  straightening law and the {AS}-{C}ohen-{M}acaulay property for quantum
  determinantal rings and quantum {G}rassmannians}, J. Algebra \textbf{301}
  (2006), no.~2, 670--702.

\bibitem{LenaganRussell}
Tom~H. Lenagan and Ewan~J. Russell, \emph{Cyclic orders on the quantum
  grassmannian}, Arab. J. Sci. Eng. Sect. C Theme Issues \textbf{33} (2008),
  no.~2, 337--350.

\bibitem{Postnikov}
Alexander Postnikov, \emph{Total positivity, {G}rassmannians, and networks},
  preprint, arXiv:math/0609764.

\bibitem{RussellThesis}
Ewan~J. Russell, \emph{Prime ideal structure of quantum grassmannians}, Ph.D.
  thesis, Edinburgh, 2008.

\bibitem{Scott-QMinors}
Joshua~S. Scott, \emph{Quasi-commuting families of quantum minors}, J. Algebra
  \textbf{290} (2005), no.~1, 204--220.

\bibitem{Williams}
Lauren~K. Williams, \emph{Enumeration of totally positive {G}rassmann cells},
  Adv. Math. \textbf{190} (2005), no.~2, 319--342.

\bibitem{Yakimov-Cyclicity}
Milen Yakimov, \emph{Cyclicity of {L}usztig's stratification of {G}rassmannians
  and {P}oisson geometry}, Noncommutative structures in mathematics and
  physics, K. Vlaam. Acad. Belgie Wet. Kunsten (KVAB), Brussels, 2010,
  pp.~259--263.

\end{thebibliography}

\normalsize

\end{document}